\documentclass{amsart}

\usepackage{amsmath}
\usepackage{mathtools}
\usepackage{amsthm,amssymb,color,comment}
\usepackage{bbm}
\usepackage{hyperref}
\usepackage[TS1,T1]{fontenc}
\usepackage{inputenc}
\usepackage{dsfont}
\usepackage{tikz}
\usepackage{enumitem}
\usepackage[noadjust]{cite}
\numberwithin{equation}{section}
\numberwithin{subsection}{section}

\theoremstyle{plain}
\newtheorem{thm}{Theorem}[section]
\newtheorem{cor}[thm]{Corollary}
\newtheorem{lem}[thm]{Lemma}
\newtheorem{prop}[thm]{Proposition}

\theoremstyle{definition}

\theoremstyle{remark}
\newtheorem{rem}[thm]{Remark}

\newcommand{\norm}[1]{\left\Vert#1\right\Vert}
\newcommand{\abs}[1]{\left\vert#1\right\vert}
\newcommand{\set}[1]{\left\{#1\right\}}
\newcommand{\R}{\mathbb R}

\newcommand{\eps}{\varepsilon}

\newcommand{\diverg}{\operatorname{div}}

\newcommand{\curl}{\operatorname{curl}}

\title[Strong Convergence of Vorticity for viscosity-dependent Navier friction]{2D Navier--Stokes with viscosity-dependent Navier slip: strong convergence for unbounded initial velocity}

\author{Josef Demmel}
\address{Department of Mathematics, Friedrich-Alexander-Universit\"at Erlangen-N\"urnberg, Cauerstr.~11, 91058 Erlangen, Germany}
\email{josef.demmel@fau.de}
\date{\today}


\begin{document}

\begin{abstract}
    For the two-dimensional incompressible Navier--Stokes equations with Navier friction on a smooth, bounded, simply connected domain, we establish strong vorticity convergence along a subsequence in the vanishing viscosity limit for initial vorticity in $L^p$ with $p>1$. The friction coefficients may grow as viscosity vanishes under an explicit growth restriction depending on $p$, which includes regimes of vanishing slip length. The main challenge, even for viscosity-independent friction, is to control the vorticity generated at the boundary. We split the vorticity into contributions from the initial data and the boundary. By combining renormalization and estimates for the normal velocity near the boundary, we show that the boundary contribution vanishes in the inviscid limit. As a consequence, the limiting Euler solutions conserve kinetic energy.
\end{abstract}

\maketitle

\section{Introduction}
Let $\Omega \subset \R^{2}$ be a bounded, simply connected domain with smooth boundary $\partial\Omega$. We consider the two-dimensional, incompressible Navier--Stokes system with Navier (friction) boundary  conditions
\begin{equation}\label{navier-stokes}
   \left\{
	\begin{aligned}
		&\partial_{t}u^{\nu} + (u^{\nu}\cdot\nabla)u^{\nu} + \nabla p^{\nu} = \nu\Delta u^{\nu}\text{ in }(0,T)\times\Omega, \\
		&\diverg u^{\nu} = 0\text{ in }(0,T)\times\Omega,\\
        &u^{\nu}\cdot n = 0\text{ on }(0,T)\times\partial\Omega,\\
        &2(Du^{\nu})_{S}n\cdot\tau + \alpha^\nu u^{\nu}\cdot \tau = 0 \text{ on }(0,T)\times\partial\Omega,\\
        &u^\nu(0,\cdot)=u_0 \text{ in }\Omega.
	\end{aligned}\right.
\end{equation}
The unknowns are the fluid velocity $u^{\nu}$ and the pressure $p^{\nu}$. We denote by $\nu>0$ the viscosity, by $n$ the unit normal vector field, by $\tau$ the tangent vector field, and by $(Dv)_{S}=\frac{1}{2}\left(\nabla v+\nabla v^{T}\right)$ the rate-of-strain tensor. The friction coefficient $\alpha^\nu\in C^2(\partial\Omega)$ is given and may depend on viscosity. \\
A key feature of the Navier boundary condition is that it implies a boundary condition for the vorticity \cite[Lemma 2.1]{clopeau} due to the formal equality
\begin{equation*}
    \curl u^\nu=2 (Du^\nu)_Sn\cdot\tau + 2\kappa u^\nu\cdot\tau \text{ on }\partial\Omega,
\end{equation*}
where $\kappa$ is the curvature of $\partial\Omega$. Applying the curl to the momentum equation, we get the vorticity formulation
\begin{equation}\label{navier-stokes-vorticity}
    \left\{
	\begin{aligned}
		&\partial_{t}\omega^{\nu} + u^{\nu}\cdot\nabla \omega^{\nu} = \nu\Delta \omega^{\nu}\text{ in }(0,T)\times\Omega, \\
		&\omega^{\nu}=\curl u^{\nu}\text{ in }(0,T)\times\Omega,\\
        &u^{\nu}\cdot n = 0\text{ on }(0,T)\times\partial\Omega,\\
        &\omega^{\nu}=(2\kappa-\alpha^\nu)u^{\nu}\cdot\tau\text{ on }(0,T)\times\partial\Omega,\\
        &\omega^\nu(0,\cdot)=\omega_0 \text{ in }\Omega.
	\end{aligned}\right.
\end{equation}
The vorticity formulation is formally equivalent, as the velocity $u^{\nu}$ can be recovered from the vorticity $\omega^{\nu}$ using the Biot--Savart law (see also Lemma \ref{curl_estimate}). Explicitly, we write
\begin{equation*}
    u^{\nu}(t)=K_{\Omega}(\omega^{\nu}(t))\text{, with }K_{\Omega}(f)(y):=\int_{\Omega}\nabla^\perp_y G_{\Omega}(x,y) f(x)\;dx,
\end{equation*}
where $G_{\Omega}$ is the Dirichlet Green's function. In particular, an initial vorticity $\omega_0\in L^p(\Omega)$ with $p>1$ identifies an initial velocity $u_0=K_{\Omega}(\omega_0)\in W^{1,p}(\Omega)\subset L^2(\Omega)$.\\
Our goal is to study the vanishing viscosity limit, i.e. $\nu\to0$, in which the Navier--Stokes system formally reduces to the Euler system
\begin{equation}\label{euler}
    \left\{
	\begin{aligned}
		&\partial_{t}u + u\cdot\nabla u + \nabla p = 0 \text{ in }(0,T)\times\Omega,\\
		&\diverg u=0 \text{ in }(0,T)\times\Omega,\\
        &u\cdot n=0 \text{ on }(0,T)\times\partial\Omega,\\
        &u(0,\cdot)=u_0 \text{ in }\Omega.
	\end{aligned}\right.
\end{equation}
The transition from Navier--Stokes to Euler as viscosity vanishes has been a classical problem with extensive research for decades (see \cite{maekawa} for an overview). Without physical boundaries, convergence in two dimensions is known for regular Euler solutions \cite{mcgrath}, and also for weak Euler solution classes with some control on the initial vorticity \cite{chemin, ciampa, seis}. With physical boundaries, the Navier--Stokes and Euler boundary conditions generally do not match, because the viscous term disappears in the limit. This mismatch leads to boundary layer effects and makes the inviscid limit more difficult. For the no-slip condition $u^\nu=0$, they are so challenging that convergence (or nonconvergence) of velocity is generally still open even in two dimensions when the limiting Euler solution is smooth. One difficulty here is that the no-slip condition does not prescribe boundary values for the vorticity, which blocks the usual compactness approach based on viscosity-independent a priori estimates for the vorticity.\\ 
As already mentioned, the Navier boundary condition, originally introduced by C. Navier in 1827 \cite{navier}, admits a formulation in terms of vorticity that makes the analysis more manageable. It can also be seen as a relaxation of the no-slip condition, as $u^\nu\cdot\tau=0$ is formally recovered for $\alpha^\nu=\infty$. For friction coefficients independent of viscosity, there are the prominent results from Clopeau, Mikeli\'{c}, and Robert \cite{clopeau}, showing strong convergence of velocity for bounded initial vorticity, and Lopes Filho, Nussenzveig Lopes, and Planas \cite{nussenzveig}, extending it to initial vorticity in $L^p(\Omega)$ for $p>2$. The proof in \cite{nussenzveig} is based on a uniform vorticity bound 
\begin{equation*}
    \norm{\omega^\nu}_{L^\infty(0,T;L^p(\Omega))}\leq C\left(\norm{\omega_0}_{L^p(\Omega)}\right) \text{ for } p>2,
\end{equation*}
where the constant $C>0$ is independent of viscosity $\nu$. The restriction $p>2$ is used to bound the boundary velocity in $L^\infty$ by the $L^p$-norm of the vorticity, an estimate that is unavailable for $p\in(1,2]$. Whether the bound can be extended to $p\in(1,2]$ is explicitly raised in \cite[Section 4, Remark 1]{nussenzveig} and was the starting point for this work, and we can remove this restriction (Corollary \ref{navier-parabolic-smooth-q=0}). Moreover, we can strengthen the velocity convergence in \cite{nussenzveig} to strong convergence of vorticity: Let $\omega_0\in L^p(\Omega)$ with $p\in(1,\infty)$. Then there exists a sequence $\nu_k\to0$ such that
\begin{equation*}
    \omega^{\nu_k}\to \curl u \text{ strongly in }C([0,T];L^p(\Omega)).
\end{equation*}
For $p>2$, strong convergence of vorticity has already been established in \cite{chemetov, demmel}.\\
When the friction coefficient is positive, $1/\alpha^\nu$ denotes the slip length. Experimental and numerical studies show that a decrease in viscosity can lead to a decrease in slip length \cite{herrero, pelz}. This motivates the study of coefficients $\alpha^\nu$ that may grow with vanishing viscosity $\nu$. We will impose the following growth bound
\begin{equation}\label{alpha-bound}
    \norm{\alpha^\nu}_{L^\infty(\partial\Omega)}\leq C\nu^{-q},\quad q\geq0,
\end{equation}
where the constant $C>0$ is independent of $\nu$.  By our choice of \eqref{alpha-bound}, we include for $q>0$ the possibility of vanishing slip (for positive friction) in the vanishing viscosity limit, which allows the boundary vorticity to grow as $\nu\to 0$. Nevertheless, our main result establishes strong vorticity convergence below a threshold on $q$ determined by the integrability of the initial vorticity.
\newpage
\begin{thm}\label{main-result}
    Let $\Omega\subset\R^2$ be a bounded, simply connected domain with smooth boundary, $T\in(0,\infty)$, $\omega_{0}\in L^{p}(\Omega)$ for $p\in(1,\infty)$, and let $u^{\nu}$ be the unique solution to \eqref{navier-stokes} for $u^{\nu}(0,\cdot)=K_{\Omega}(\omega_{0})$ with $\alpha^\nu\in C^2(\partial\Omega)$ which satisfies \eqref{alpha-bound}. Then, for $\rho\in(1,p]$ and $q\geq 0$ satisfying $\rho q<\min\set{\frac{p-1}{2p},\frac{1}{4}}$, there exists a sequence $\nu_{k}\to 0$ such that 
    \begin{equation*}
        u^{\nu_{k}}\to u \text{ strongly in } C([0,T];W^{1,\rho}(\Omega))\text{ as }k\to\infty,\text{ }
    \end{equation*}
    and $u\in C([0,T];W^{1,p}(\Omega))$ is a weak solution of \eqref{euler} for $u(0,\cdot)=K_\Omega(\omega_{0})$. In particular,
    \begin{equation*}
        \omega^{\nu_{k}}\to \curl u \text{ strongly in } C([0,T];L^{\rho}(\Omega))\text{ as }k\to\infty,
    \end{equation*}
    and the limiting vorticity $\omega=\curl u$ is a renormalized solution of $\partial_t\omega+u\cdot\nabla\omega=0$ with initial datum $\omega_0$. Consequently,
    \begin{equation*}
        \norm{\omega(t)}_{L^p(\Omega)}=\norm{\omega_0}_{L^p(\Omega)} \text{ for every }t\in[0,T].
    \end{equation*}
\end{thm}
The underlying estimates also yield conservation of kinetic energy for the limiting Euler solution. Analogous results for vanishing-viscosity limits without physical boundaries were established in \cite{cheskidov, ciampa}.
\begin{cor}\label{energy-conservation}
    Under the same assumptions as in Theorem \ref{main-result}, we have
    \begin{equation*}
        \nu\int_0^T\norm{\nabla u^\nu(t)}_{L^2(\Omega)}^2\;dt+\nu\int_0^T\int_{\partial\Omega}(\abs{\kappa}+\abs{\alpha^\nu})\abs{u^\nu(t)\cdot\tau}^2\;dS\;dt\to0 \text{ as }\nu\to0.
    \end{equation*}
    In particular, every limiting Euler solution satisfies
    \begin{equation*}
        \norm{u(t)}_{L^2(\Omega)}=\norm{K_{\Omega}(\omega_0)}_{L^2(\Omega)}\text{ for every }t\in[0,T].
    \end{equation*}
\end{cor}
For viscosity-dependent Navier friction, strong convergence of velocity towards sufficiently regular Euler solutions has been established in \cite{kim1, kim2, paddick} by classical energy estimates. For analytic data in the half-plane, short-time convergence of velocity at the critical scale $q=1$ was shown in \cite{nguyen}.  Also, Kato-type convergence criteria have been developed in \cite{maekawa2}. Regarding weak Euler solutions, a conditional convergence result for velocity was shown in \cite{drivas}, assuming viscosity-independent interior bounds on the viscous flows.\\
For initial vorticity in $L^p(\Omega)$ with $p>1$, we establish unconditional strong convergence of vorticity to the vorticity of a weak Euler solution. Whenever $q<\min\set{\frac{p-1}{2p},\frac{1}{4}}$, the vorticities converge along a subsequence in $C([0,T];L^\rho(\Omega))$ for some $\rho>1$ to the vorticity of a weak solution to Euler for the same initial data. This convergence can fail already at $q=\frac{1}{2}$ (Lemma \ref{counterexample}). The counterexample uses the same boundary obstruction as in Kelliher's analysis of the no-slip problem \cite{kelliher2}: strong vorticity convergence implies strong convergence of the tangential velocities at the boundary, which fails in our counterexample. The restriction on $q$ in Theorem \ref{main-result} originates from the methods and is not expected to be optimal. We now explain how we control the vorticity generated at the boundary.
\subsection{Proof strategy}
A main theme in our argument is to divide the vorticity $\omega^\nu$ into an initial contribution $a^\nu$, which solves
\begin{equation}\label{initial-vorticity}
        \left\{
	\begin{aligned}
		&\partial_{t}a^{\nu} + u^{\nu}\cdot\nabla a^{\nu} = \nu\Delta a^{\nu}\text{ in }(0,T)\times\Omega, \\
        &a^{\nu}=0\text{ on }(0,T)\times\partial\Omega,\;a^\nu(0)=\omega_0 \text{ in }\Omega,
	\end{aligned}\right.
    \end{equation}
and a boundary contribution $b^\nu$, which solves
\begin{equation}\label{boundary-vorticity}
        \left\{
	\begin{aligned}
		&\partial_{t}b^{\nu} + u^{\nu}\cdot\nabla b^{\nu} = \nu\Delta b^{\nu}\text{ in }(0,T)\times\Omega, \\
        &b^{\nu}=(2\kappa-\alpha^\nu)u^\nu\cdot\tau\text{ on }(0,T)\times\partial\Omega,\;b^\nu(0)=0 \text{ in }\Omega.
	\end{aligned}\right.
\end{equation}
With $u^\nu$ regarded as a prescribed drift for both auxiliary problems, linearity and uniqueness recover $\omega^\nu=a^\nu+b^\nu$. The homogeneous boundary condition allows us to control $a^\nu$ by contractivity (Lemma \ref{lp-contraction}) and parabolic smoothing (Lemma \ref{parabolic-smoothing}). Therefore, the main challenge lies in controlling $b^\nu$. While the boundary values of $b^\nu$ might grow in the vanishing viscosity limit, we will show that their effect inside the domain still vanishes. \\
First, we establish an estimate combining renormalization with control of the boundary for advection-diffusion equations with divergence-free drift tangent to the boundary (Lemma \ref{advection-diffusion-lemma}). We assume a suitable growth bound on the inward normal component of the drift in a boundary collar \eqref{normal-drift-assumption}, which lets us replace the uncontrolled diffusive boundary flux (see also \cite[Section 4, Remark 2]{nussenzveig}) in the renormalization identity by a manageable boundary contribution. Second, we will show that $u^\nu$ satisfies \eqref{normal-drift-assumption}, in a sufficiently small boundary collar, by the bound
\begin{equation}\label{normal-drift-bound-intro}
    \abs{u^\nu(s)\cdot \nabla d}\leq C(\nu s)^{-1/r}\norm{\curl u_0}_{L^r(\Omega)}d^\gamma,\quad s\in(0,T),
\end{equation}
where $r\in (1,\min\set{2,p})$, $\gamma\in(0,1)$ will be chosen below, and $d(x)=\text{dist}(x,\partial\Omega)$ is the distance function to the boundary. To obtain the normal drift bound \eqref{normal-drift-bound-intro}, we first establish parabolic smoothing estimates for the vorticity (Lemma \ref{navier-parabolic-smoothing}). Then, we use the Biot--Savart law and Morrey's inequality to deduce Hölder estimates for the velocity, which imply the desired estimate thanks to the impermeability condition (Corollary \ref{normal-velocity-bound-cor}). The restriction $r<2$ arises in Lemma \ref{navier-parabolic-smoothing} from combining the interpolation estimate (Lemma \ref{biot-savart-tangential-estimates}) with a uniform $L^1(\Omega)$ bound for the newly generated vorticity, so this part of the argument does not benefit from higher initial integrability.\\
With these tools, we obtain a uniform vorticity bound (Proposition \ref{compatible-uniform-vorticity-bound} and Lemma \ref{vorticity-bound-weak}), and consequently, for $\rho\in(1,p]$,
\begin{equation*}
    \sup_{t\in[0,T]}\norm{b^\nu(t)}_{L^\rho(\Omega)}^\rho
\le C\norm{\omega_0}_{L^\rho(\Omega)}^\rho
\left(
\nu^{1/2-\rho q}
+
\nu^{(\gamma-1/r)/(1+\gamma)-\rho q}
\norm{\omega_0}_{L^r(\Omega)}^{1/(1+\gamma)}
\right).
\end{equation*}
The first summand comes from diffusion and requires $\rho q<1/2$ to vanish. Thus, this term can decay for some $\rho>1$ only when $q<1/2$, consistent with the failure of strong vorticity convergence in our counterexample for $q=1/2$. The second summand comes from the normal drift bound \eqref{normal-drift-bound-intro} and imposes the stronger restriction
\begin{equation*}
    \rho q<\sup_{\substack{\gamma\in(0,1)\\r\in(1,\min\set{2,p})}}\frac{\gamma-1/r}{1+\gamma}=\min\set{\frac{p-1}{2p},\frac{1}{4}}.
\end{equation*}
Consequently, $b^\nu$ vanishes uniformly in time in $L^\rho(\Omega)$. Finally, compactness gives velocity convergence along a subsequence, and stability of the homogeneous problem (Lemma \ref{linear-stability}) yields strong convergence of $a^\nu$ to the limiting Euler vorticity. Together, this proves the asserted strong vorticity convergence.
\subsection{Declaration of AI Use}
The author directed the mathematical work on this project. Large language models, including OpenAI’s GPT-5.5 Pro and GPT-5.6 Sol, served as research assistants, advising the author on specific questions that arose in developing the proof. The counterexample in Lemma \ref{counterexample} was found by OpenAI's GPT-6 Astra.\\
The author wrote the manuscript from scratch and takes full responsibility for the mathematical content and presentation of the paper.
\section{Notation}
In the following, $x,y\in\overline{\Omega}$ will be spatial variables, while $t,s\in[0,T]$ are reserved for time variables. We denote with $v\cdot w$ the Euclidean scalar product and with $A:B$ the scalar product of matrices. In two dimensions, we set
\begin{equation*}
    \curl v := \partial_{x_1}v_2-\partial_{x_2}v_1,\quad\nabla^\perp f := \binom{-\partial_{x_2}}{\partial_{x_1}}f.
\end{equation*}
Regarding the boundary, we set $\tau=n^\perp=\binom{-n_2}{n_1}$ and the curvature $\kappa$ by $\partial_\tau n=\kappa\tau$. We write $\partial_nf=\nabla f\cdot n$ for the normal derivative, $d(x)=\text{dist}(x,\partial\Omega)$ for the distance to the boundary, and $\Omega_\xi:=\set{x\in\Omega:\text{dist}(x,\partial\Omega)<\xi}$ for a boundary collar of size $\xi>0$.
We normalize the Dirichlet Green function by $\Delta_yG_\Omega(x,\cdot)=\delta_x$ and $G_\Omega(x,\cdot)=0$ on $\partial\Omega$.
Following \cite{clopeau}, we introduce the following Hilbert spaces
\begin{equation*}
    \begin{aligned}
        &H=\set{v\in L^{2}(\Omega): \diverg v=0 \text{ in } \Omega \text{ and }v\cdot n = 0 \text{ on }\partial\Omega},\\
        &V=\set{v\in H^{1}(\Omega): \diverg v=0 \text{ in } \Omega \text{ and }v\cdot n = 0 \text{ on }\partial\Omega},\\
        &W=\set{v\in H^{2}(\Omega)\cap V: 2(Dv)_{S}n\cdot\tau + \alpha^\nu v\cdot\tau=0 \text{ on }\partial\Omega}.
    \end{aligned}
\end{equation*}
For $W$, we will suppress the dependence on $\alpha^\nu$ in the notation. Note that the boundary conditions in these definitions are intended in an appropriate trace sense. For $\gamma\in(0,1)$, the Hölder seminorm is
\begin{equation*}
    [v]_{C^0,\gamma(\overline{\Omega})}:= \sup_{\substack{x,y\in\overline{\Omega}\\ x\neq y}}\frac{\abs{v(x)-v(y)}}{\abs{x-y}^\gamma}.
\end{equation*}
We will use the same notation for scalar- and vector-valued function spaces. Finally, values of constants in an inequality may change from line to line and may depend on the domain and also on functions of the domain.
\section{Weak and strong solutions for Navier friction}
For the convenience of the reader, we repeat well-posedness and regularity results for Navier boundary conditions. Like in the Dirichlet case, there exists a unique Leray--Hopf solution for $u_0\in H$ that satisfies the energy inequality.
\begin{thm}[Theorem 6.1, \cite{kelliher}]\label{energybound}
    For $\nu>0$ and $u_{0}\in H$, there exists a unique weak solution $u^\nu$ of \eqref{navier-stokes}, i.e. a vector field 
    \begin{equation*}
        u^\nu\in C([0,T];H)\cap L^2(0,T;V)
    \end{equation*}
    with $u^\nu(0,\cdot)=u_0$, satisfying
    \begin{equation}\label{weak-formulation}
        \frac{d}{dt}\int_{\Omega}u^\nu\cdot v dx+\int_{\Omega}(u^\nu\cdot\nabla u^\nu)\cdot vdx+\nu\int_{\Omega}\nabla u^\nu:\nabla vdx=\nu\int_{\partial\Omega}(\kappa-\alpha^\nu)u^\nu\cdot v dS,
    \end{equation}
    for every $v\in V$. Moreover, $\partial_t u^\nu\in L^2(0,T;V')$, and for every $t\in[0,T]$, the following energy identity holds
    \begin{equation}\label{energy-identity}
        \norm{u^\nu(t)}_{L^2(\Omega)}^2+2\nu\int_0^t\norm{\nabla u^\nu}_{L^2(\Omega)}^2\;ds= \norm{u_0}_{L^2(\Omega)}^2+2\nu\int_0^t\int_{\partial\Omega}(\kappa-\alpha^\nu)\abs{u^\nu\cdot\tau}^2\;dS\;ds.
    \end{equation}
\end{thm}
\begin{rem}
    For $\alpha^\nu\geq 0$, we have the stronger energy inequality
    \begin{equation*}
        \norm{u^\nu(t)}_{L^2(\Omega)}^2+\nu\int_0^t\norm{\nabla u^\nu}_{L^2(\Omega)}^2\;ds\leq e^{C\nu t}\norm{u_0}_{L^2(\Omega)}^2.
    \end{equation*}
    For positive times, the weak solution of Theorem \ref{energybound} satisfies
    \begin{equation*}
        u^\nu\in C([\eps,T];V)\cap L^2(\eps,T;W) \text{ for every }\eps\in(0,T).
    \end{equation*}
    Since $u^\nu\in L^2(0,T;V)$, we may choose $s\in(0,\eps)$ with $u^\nu(s)\in V$ and apply \cite[Theorem 1.1]{demmel} starting at time $s$. 
\end{rem}

A common idea, which refers back to \cite{nussenzveig}, is to consider approximations of the initial data that satisfy the Navier boundary condition in order to work with strong solutions.
\begin{lem}[\cite{nussenzveig}]\label{compatible-lemma}
    Let $\omega_0\in L^p(\Omega)$ for $p>1$. Then there exists a sequence of compatible functions $\set{\omega_{0,n}}$, i.e., $\{\omega_{0,n}\}\subset H^{1}(\Omega)\cap L^{\infty}(\Omega)$ with $K_{\Omega}(\omega_{0,n})\in W$, that converges strongly to $\omega_0$ in $L^p(\Omega)$.
\end{lem}

\begin{thm}[Theorem 2.3, \cite{clopeau}; Theorem 7.1, \cite{kelliher}]\label{theorem-clopeau}
    For $\nu>0$ and $\omega_{0}\in H^1(\Omega)\cap L^\infty(\Omega)$ with $u_{0}=K_{\Omega}(\omega_{0})\in W$, there exists a unique strong solution $u^\nu\in C([0,T];W)$ to \eqref{navier-stokes} with $u^\nu(0,\cdot)=K_\Omega(\omega_0)$. Moreover, $\partial_{t}u^\nu\in L^{2}(0,T;V)\cap C([0,T];H)$ and the associated vorticity $\omega^\nu=\curl u^\nu$ is bounded, $\omega^\nu\in C([0,T];H^{1}(\Omega))\cap L^{\infty}((0,T)\times\Omega)$.
\end{thm}

\section{A priori estimates for the advection-diffusion equation}
If we split the vorticity $\omega^\nu$ into an initial contribution $a^\nu$ and a boundary contribution $b^\nu$, they individually are no longer solutions to \eqref{navier-stokes-vorticity}, only to the more general advection-diffusion equation below.
\begin{equation}\label{advection-diffusion equation}
    \left\{
	\begin{aligned}
		&\partial_{t}\theta^\nu + v^\nu\cdot\nabla \theta^\nu = \nu\Delta \theta^\nu\text{ in }(0,T)\times\Omega, \\
		&\diverg v^\nu =0\text{ in }(0,T)\times\Omega,\\
        &v^\nu\cdot n = 0\text{ on }(0,T)\times\partial\Omega,\\
        &\theta^\nu=g^\nu\text{ on }(0,T)\times\partial\Omega.
	\end{aligned}\right.
\end{equation}
In this section, we quantify the effect of the boundary data $g^\nu$ under a suitable assumption on the drift $v^\nu$ (Lemma \ref{advection-diffusion-lemma}) and establish well-known a priori estimates in the case $g^\nu=0$ (Lemma \ref{lp-contraction} and Lemma \ref{parabolic-smoothing}). For a simpler notation and as all statements are for fixed $\nu\in(0,1)$, we will omit the superscript $\nu$ throughout this section.
\begin{lem}\label{advection-diffusion-lemma}
    Fix $\nu\in(0,1)$, $\gamma>0$, and a convex function $\beta\in C(\R;[0,\infty))$. For 
    $v\in L^2_{loc}(0,T;H)$ and $g\in L^2_{loc}(0,T;H^{1/2}(\partial\Omega))$ with $\beta(g)\in L^\infty(0,T;L^1(\partial\Omega))$, suppose there exists a corresponding weak solution $\theta$ of \eqref{advection-diffusion equation} satisfying
    \begin{equation*}
        \theta\in C([0,T];L^1(\Omega))\cap L^2_{loc}(0,T;H^1(\Omega)) \text{ with }\beta(\theta(0))\in L^1(\Omega).
    \end{equation*} Assume also that there exists a strictly positive, nonincreasing $B\in W^{1,1}_{loc}(0,T)$ with $B^{1/(1+\gamma)}\in L^1(0,T)$ such that
    \begin{equation}\label{normal-drift-assumption}
        v(s,x)\cdot\nabla d(x)\leq B(s)d(x)^\gamma \text{ for almost every }(s,x)\in(0,T)\times \Omega_{\xi},
    \end{equation}
    where the collar width $\xi>0$ is fixed independently of $v$ and $\nu$ and is small enough that $d\in C^2(\Omega_{\xi})$ and $\xi\norm{\Delta d}_{L^\infty(\set{d<\xi})}\leq \frac{1}{4}.$ Then, for every $0\leq s<t\leq T$, we have
    \begin{equation*}
        \int_\Omega \beta(\theta(t))\;dx\leq \int_\Omega \beta(\theta(s))\;dx+C\int_{s}^t\left(\sqrt{\frac{\nu}{t-r}}+\nu^{\gamma/(1+\gamma)}B(r)^{1/(1+\gamma)}\right) \int_{\partial\Omega} \beta(g(r))\;dS\;dr,
    \end{equation*}
    where the constant $C>0$ is independent of the drift $v$ and the viscosity $\nu$.
\end{lem}
\begin{proof}
    The idea is to test the advection-diffusion equation of \eqref{advection-diffusion equation} with $\beta'(\theta)\varphi$, using \eqref{normal-drift-assumption} to choose $\varphi$ so that diffusion can absorb the other bulk terms at the cost of a manageable boundary contribution. Fix a convex $\beta\in C^2(\R;[0,\infty))$ with $\beta',\beta''\in L^{\infty}(\R)$. In the end, the general case will follow by approximation. To justify the chain rule for the time derivative with $\beta'(\theta)\varphi$, we mollify in time. Fix $0\leq s<t\leq T$ and choose $s<a<b<t$. Let $\rho\in C^\infty_c((-1,1))$ be nonnegative with $\int_\R\rho=1$. For $\varepsilon\in(0,\min\set{a,T-b})$ and $r\in[a,b]$, we set
    \begin{equation*}
        f_\varepsilon(r):=\int_0^T\rho_\varepsilon(r-\tau)f(\tau)\;d\tau, \text{ with }\rho_\varepsilon(r):=\frac{1}{\varepsilon}\rho\left(\frac{r}{\varepsilon}\right).
    \end{equation*}
    The mollified $\theta_\varepsilon$ solves
    \begin{equation}\label{mollified-advection-diffusion}
        \partial_r\theta_\varepsilon
        + (v\cdot\nabla\theta)_\varepsilon=\nu\Delta\theta_\varepsilon.
    \end{equation}
    For nonnegative $\varphi\in W^{1,1}(a,b;L^\infty(\Omega))\cap L^2(a,b;H^1_0(\Omega))$ with $\Delta \varphi\in L^\infty((a,b)\times\Omega)$, we can test \eqref{mollified-advection-diffusion} with $\beta'(\theta_\varepsilon)\varphi$, and use integration by parts to get
    \begin{align*}
        &\int_\Omega \beta(\theta_\varepsilon(b))\varphi(b)\;dx-\int_\Omega \beta(\theta_\varepsilon(a))\varphi(a)\;dx-\int_a^b\int_{\Omega}\beta(\theta_\varepsilon)\partial_r\varphi\;dx\;dr+\int_a^b\int_{\Omega}(v\cdot\nabla\theta)_\varepsilon\beta'(\theta_\varepsilon)\varphi\;dx\;dr\\
        &=-\nu\int_a^b\int_{\Omega}\nabla\beta(\theta_\varepsilon)\cdot\nabla\varphi\;dx\;dr-\nu\int_a^b\int_{\Omega}\abs{\nabla\theta_\varepsilon}^2\beta''(\theta_\varepsilon)\varphi\;dx\;dr.
    \end{align*}
    From the convexity of $\beta$, we know that $\beta''\geq0$, and thus the last term on the right-hand side is nonpositive. The standard properties of time mollification, together with the assumed $L^1(\Omega)$-continuity of $\theta$, yield
    \begin{align*}
        \theta_\eps\to\theta
        \text{ in }L^2(a,b;H^1(\Omega))\cap C([a,b];L^1(\Omega))\text{ and }
        (v\cdot\nabla\theta)_\eps\to v\cdot\nabla\theta
        \text{ in }L^1((a,b)\times\Omega).
    \end{align*}
    For a subsequence, we also have $\theta_\eps\to\theta$ almost everywhere. Therefore, by the choice of $\beta$, we can pass to the limit in all the terms to obtain
    \begin{align*}
        &\int_\Omega \beta(\theta(b))\varphi(b)\;dx-\int_\Omega \beta(\theta(a))\varphi(a)\;dx-\int_a^b\int_{\Omega}\beta(\theta)\partial_r\varphi\;dx\;dr+\int_a^b\int_{\Omega}(v\cdot\nabla\theta)\beta'(\theta)\varphi\;dx\;dr\\
        &\leq-\nu\int_a^b\int_{\Omega}\nabla\beta(\theta)\cdot\nabla\varphi\;dx\;dr.
    \end{align*}
    With integration by parts and the Dirichlet boundary condition of \eqref{advection-diffusion equation}, this yields
    \begin{align*}
        \int_\Omega \beta(\theta(b))\varphi(b)\;dx\leq&\int_\Omega \beta(\theta(a))\varphi(a)\;dx+\int_a^b\int_{\Omega}\beta(\theta)[\partial_r\varphi+v\cdot\nabla\varphi+\nu\Delta\varphi]\;dx\;dr\\&-\nu\int_a^b\int_{\partial\Omega}\beta(g)\partial_n\varphi\;dS\;dr.
    \end{align*}
    Now, we want to construct $\varphi$. We seek a nonnegative $\varphi \in W^{1,1}(a,b;L^\infty(\Omega))\cap L^2(a,b;H^1_0(\Omega))$ with $\Delta \varphi\in L^\infty((a,b)\times\Omega)$ satisfying the backward adjoint inequality, i.e. $\partial_r\varphi+v\cdot\nabla\varphi+\nu\Delta\varphi\leq 0$, being bounded by 1 in $(s,t)\times\Omega$, converging pointwise in $\Omega$ to 1  as $r\to t$, and with controllable boundary flux $-\nu\partial_n\varphi$. An appropriate choice that meets the regularity criteria is
    \begin{equation*}
        \varphi(r,x):=\begin{cases}
            2\frac{d(x)}{l(r)}-\frac{d(x)^2}{l(r)^2},&d(x)<l(r),\\
            1,&d(x)\geq l(r).
        \end{cases}\text{ with }\varphi=0 \text{ and}-\nu\partial_n\varphi=\frac{2\nu}{l(r)} \text{ on }\partial\Omega,
    \end{equation*}
    where
    \begin{equation*}
        l(r):=\min\set{\xi,\sqrt{\nu(t-r)},\left(\frac{\nu}{4B(r)}\right)^{1/(1+\gamma)}}.
    \end{equation*}
    Let us explain the choice of $l(r)$. We need $l(r)\leq \xi$ to ensure that the distance function $d$ is smooth and the curvature contribution $\nu\Delta d$ is controllable. The second condition, $l(r)\leq \sqrt{\nu(t-r)}$, has two purposes. First, it ensures the convergence $\lim_{r\to t}\varphi(r,x)=1$ for $x\in\Omega$. Second, it yields the bound $-l'(r)\leq\frac{\nu}{2l(r)}$, as the other two terms are nondecreasing. To understand the last term, it is best to look at the backward adjoint problem for $d(x)<l(r)$, where we also use \eqref{normal-drift-assumption}:
    \begin{align*}
        \partial_r\varphi+v\cdot\nabla\varphi+\nu\Delta\varphi&=\frac{2}{l}\underbrace{\left(1-\frac{d}{l}\right)}_{\in(0,1)}\left(\underbrace{-\frac{d}{l}l'}_{\leq\nu/(2l)}+\underbrace{v\cdot\nabla d}_{\leq Bd^\gamma}+\underbrace{\nu\Delta d}_{\leq\nu/(4l)}\right)-\frac{2\nu}{l^2}\\
        &\leq\frac{2}{l}B d^{\gamma}-\frac{\nu}{2l^2}\leq0,
    \end{align*}
    due to the third term
    \begin{equation*}
        \frac{2}{l}B d^{\gamma}\leq \frac{\nu}{2l^2}\Leftrightarrow \underbrace{l d^\gamma}_{\leq l^{\gamma+1}}\leq \frac{\nu}{4B}.
    \end{equation*}
    For $d(x)\geq l(r)$, the constant $1$ solves the backward adjoint problem exactly, and thus we have
    \begin{align*}
        &\int_\Omega \beta(\theta(b))\varphi(b)\;dx\leq\int_\Omega \beta(\theta(a))\underbrace{\varphi(a)}_{\leq 1}\;dx+\nu\int_a^b\frac{2}{l}\int_{\partial\Omega}\beta(g)\;dS\;dr\\
        &\leq\int_\Omega \beta(\theta(a))\;dx+2\nu\int_a^b\left(\frac{1}{\xi}+\frac{1}{\sqrt{\nu(t-r)}}+\left(\frac{4B(r)}{\nu}\right)^{1/(1+\gamma)}\right)\int_{\partial\Omega}\beta(g)\;dS\;dr.\\
        &\leq\int_\Omega \beta(\theta(a))\;dx+C\int_{s}^t\left(\sqrt{\frac{\nu}{t-r}}+\nu^{\gamma/(1+\gamma)}B(r)^{1/(1+\gamma)}\right)\int_{\partial\Omega}\beta(g)\;dS\;dr.
    \end{align*}
    The Lipschitz continuity of $\beta$ gives $\beta(\theta)\in C([s,t];L^1(\Omega))$, and therefore we can recover the endpoints
    \begin{equation*}
        \lim_{a\to s}\int_\Omega \beta(\theta(a))\;dx=\int_\Omega \beta(\theta(s))\;dx \text{ and }\lim_{b\to t}\int_\Omega \beta(\theta(b))\varphi(b)\;dx=\int_\Omega \beta(\theta(t))\;dx.
    \end{equation*}
    In total, we have shown
    \begin{equation}\label{c2-beta-inequality}
        \int_\Omega \beta(\theta(t))\;dx\leq \int_\Omega \beta(\theta(s))\;dx+C\int_{s}^t\left(\sqrt{\frac{\nu}{t-r}}+\nu^{\gamma/(1+\gamma)}B(r)^{1/(1+\gamma)}\right) \int_{\partial\Omega} \beta(g)\;dS\;dr,
    \end{equation}
    so it remains to remove the extra smoothness assumptions on $\beta$. For general convex $\beta\in C(\R;[0,\infty))$, we approximate $\beta$ by convex $\beta_k\in C^2(\R;[0,\infty))$ with $\beta'_k,\beta_k''\in L^{\infty}(\R)$, $0\leq\beta_k\leq\beta+\frac{1}{k}$, and $\beta_k\to\beta$ pointwise. First use \eqref{c2-beta-inequality} with $\beta_k$ and $s=0$, because we only know $\beta(\theta(0))\in L^1(\Omega)$. Fatou's lemma on the left-hand side and dominated convergence on the right-hand side then yield the desired estimate starting from $0$. As the right-hand side is finite, we get $\beta(\theta(t))\in L^1(\Omega)$ for every $t\in[0,T]$. Now, we can pass to the limit for arbitrary $s\in[0,t)$, which proves the assertion.
\end{proof}
For the sake of completeness, we will repeat some properties of the advection-diffusion equation with homogeneous Dirichlet data, i.e. $\theta=0$ on $(0,T)\times\partial\Omega$. We start with the contraction property for all $L^p(\Omega)$ norms. 
\begin{lem}\label{lp-contraction}
    Fix $\nu\in(0,1)$ and $p\in[1,\infty]$. Let $v\in L^2(0,T;H)$, and 
    \begin{equation*}
        \theta \in C([0,T];L^1(\Omega))\cap L^2_{loc}(0,T;H^1_0(\Omega))\text{ with } \theta(0)\in L^p(\Omega),
    \end{equation*}
     be a corresponding weak solution to \eqref{advection-diffusion equation} with $g=0$. Then, for every $0\leq s<t\leq T$, it holds that
     \begin{align*}
         &\norm{\theta(t)}_{L^p(\Omega)}\leq \norm{\theta(s)}_{L^p(\Omega)}.
     \end{align*}
\end{lem}
\begin{proof}
    Let $\beta\in C^2(\R;[0,\infty))$ be convex with $\beta',\beta''\in L^\infty(\R)$, and $\beta'(0)=0$, which ensures that $\beta'(\theta)$ vanishes on the boundary. Testing with $\beta'(\theta)$ and integration by parts, justified by the time mollification argument as in the proof of Lemma \ref{advection-diffusion-lemma}, yields
    \begin{equation}\label{0-boundary-convex-estimate}
        \frac{d}{dr}\int_\Omega\beta(\theta(r))\;dx=-\nu\int_\Omega \beta''(\theta(r))\abs{\nabla\theta(r)}^2\;dx\leq 0\text{ in }\mathcal{D}'(0,T).
    \end{equation}
    Since $\beta$ is Lipschitz, we have $\beta(\theta)\in C([0,T];L^1(\Omega))$, and thus 
    \begin{equation}\label{convex-contraction}
        \int_\Omega\beta(\theta(t))\;dx\leq \int_\Omega\beta(\theta(s))\;dx \text{ for every }0\leq s<t\leq T.
    \end{equation}
    For $p\in[1,\infty)$, choose convex functions $\beta_k\in C^2(\R;[0,\infty))$ satisfying
    \begin{equation*}
        \beta'_k(0)=0,\quad \beta'_k,\beta''_k\in L^\infty(\R),\quad0\leq\beta_k(z)\leq \abs{z}^p+\frac{1}{k},\quad \beta_k\to\abs{\cdot}^p.
    \end{equation*}
    Like before, we get with Fatou's lemma and dominated convergence that
    \begin{equation*}
        \int_{\Omega}\abs{\theta(t)}^p\;dx\leq \liminf_{k\to\infty}\int_{\Omega}\beta_k(\theta(t))\;dx\leq \lim_{k\to\infty}\int_{\Omega}\beta_k(\theta(0))\;dx=\int_{\Omega}\abs{\theta(0)}^p\;dx,
    \end{equation*}
    hence $\theta(t)\in L^p(\Omega)$ for all $t\in[0,T]$. Repeating the argument with arbitrary $s\in[0,t)$ on the right-hand side yields
    \begin{equation}\label{finite-p-contraction}
        \norm{\theta(t)}_{L^p(\Omega)}\leq\norm{\theta(s)}_{L^p(\Omega)}.
    \end{equation}
    For $p=\infty$, we have $\theta(0)\in L^q(\Omega)$ for every finite $q$. Since $\abs{\Omega}<\infty$, letting $q\to\infty$ in the $L^q$ estimate \eqref{finite-p-contraction} yields the $L^\infty$ contraction.
\end{proof}
For positive times, we even gain integrability thanks to parabolic smoothing estimates, whose underlying argument goes back to Nash \cite{nash}. On the whole space, similar parabolic smoothing estimates have also been shown in \cite{carlen, silvestre}.   
\begin{lem}\label{parabolic-smoothing}
    Fix $\nu\in(0,1)$ and $p\in[1,\infty]$. Let $v\in L^2(0,T;H)$, and 
    \begin{equation*}
        \theta \in C([0,T];L^1(\Omega))\cap L^2_{loc}(0,T;H^1_0(\Omega))
    \end{equation*}
     be a corresponding weak solution to \eqref{advection-diffusion equation} with $g=0$. Then, for every $0\leq s<t\leq T$ and $1\leq m\leq p$ such that $\theta(s)\in L^m(\Omega)$, we have $\theta(t)\in L^p(\Omega)$ and
     \begin{align*}
         \norm{\theta(t)}_{L^p(\Omega)}\leq C[\nu(t-s)]^{-(1/m-1/p)}\norm{\theta(s)}_{L^m(\Omega)}.
     \end{align*}
     Here $1/\infty=0$, and the constant $C>0$ is independent of $v$ and the viscosity $\nu$.
\end{lem}
\begin{proof}
    For $m=\infty$, the statement is implied by Lemma \ref{lp-contraction}. Let $m\in[1,\infty)$ and first suppose $\theta(s)\in L^\infty(\Omega)$. Lemma \ref{lp-contraction} and $\theta\in C([0,T];L^1(\Omega))$ imply that $\theta\in C([s,T];L^q(\Omega))$ for every finite $q\in[1,\infty)$. For $q\geq2$, applying \eqref{0-boundary-convex-estimate} with an admissible $\beta$ that coincides with $\abs{\cdot}^q$ on the range of $\theta$ gives
    \begin{equation}\label{parabolic-smoothing-identity}
        \frac{d}{dr}\norm{\theta(r)}^q_{L^q(\Omega)}+\frac{4\nu(q-1)}{q}\norm{\nabla\abs{\theta(r)}^{q/2}}^2_{L^2(\Omega)}=0 \text{ for almost every }r\in(s,T).
    \end{equation}
    Nash's inequality \cite[p. 936]{nash}, applied to the zero extension of $f\in H^1_0(\Omega)$, gives
    \begin{equation*}
        \norm{f}_{L^2(\Omega)}^4\leq C \norm{f}^2_{L^1(\Omega)}\norm{\nabla f}^2_{L^2(\Omega)}.
    \end{equation*}
    If $\theta(t)=0$, the estimate is trivial. Otherwise, $\theta(r)\neq 0$ for every $r\in[s,t]$ by Lemma \ref{lp-contraction}. Fix $s\leq a<b\leq t$. For almost every $r\in(a,b)$, applying Nash's inequality to $f=\abs{\theta(r)}^{q/2}\in H^1_0(\Omega)$ and using $L^{q/2}$ contraction gives
     \begin{equation*}
         \norm{\theta(r)}^{2q}_{L^q(\Omega)}\leq C\norm{\theta(a)}^q_{L^{q/2}(\Omega)}\norm{\nabla\abs{\theta(r)}^{q/2}}^2_{L^2(\Omega)}.
     \end{equation*}
     With \eqref{parabolic-smoothing-identity}, we obtain, for almost every $r\in(a,b)$,
     \begin{align*}
         \frac{d}{dr}\left(\norm{\theta(r)}_{L^q(\Omega)}^{-q}\right)&=\frac{-\frac{d}{dr}\norm{\theta(r)}^q_{L^q(\Omega)}}{\norm{\theta(r)}^{2q}_{L^q(\Omega)}}=\frac{4\nu(q-1)}{q}\frac{\norm{\nabla\abs{\theta(r)}^{q/2}}^2_{L^2(\Omega)}}{\norm{\theta(r)}^{2q}_{L^q(\Omega)}}\\&\geq \frac{4\nu(q-1)}{q}\frac{1}{C\norm{\theta(a)}^q_{L^{q/2}(\Omega)}}\geq \frac{2\nu}{C\norm{\theta(a)}^q_{L^{q/2}(\Omega)}}.
     \end{align*}
     Integrating from $a$ to $b$ and using continuity at the endpoints yields
     \begin{equation}\label{q-inequality}
         \norm{\theta(b)}_{L^q(\Omega)}\leq \left(\frac{C}{\nu(b-a)}\right)^{1/q}\norm{\theta(a)}_{L^{q/2}(\Omega)}.
     \end{equation}
     The constant $C$ is independent of $q$. Set $q_k=2^km$ and $t_k=s+(t-s)(1-2^{-k})$ for $k\geq0$. Iterating \eqref{q-inequality} with $q=q_j$ on $[t_{j-1},t_j]$ for $j\in\set{1,...,k}$, we derive
    \begin{align*}
        \norm{\theta(t)}_{L^{q_k}(\Omega)}\leq \norm{\theta(t_k)}_{L^{q_k}(\Omega)}&\leq \norm{\theta(s)}_{L^{m}(\Omega)}\prod_{j=1}^k\left(\frac{C2^j}{\nu(t-s)}\right)^{1/(2^jm)}.
    \end{align*}
    By letting $k\to\infty$, we conclude
    \begin{equation*}
        \norm{\theta(t)}_{L^{\infty}(\Omega)}\leq \norm{\theta(s)}_{L^{m}}\left(\frac{4C}{\nu(t-s)}\right)^{1/m},
    \end{equation*}
    where we used for the right-hand side the series $\sum_{j=1}^\infty 2^{-j}=1$ and $\sum_{j=1}^\infty j2^{-j}=2$. For finite $p$, interpolation and Lemma \ref{lp-contraction} yield the desired estimate
    \begin{align*}
        \norm{\theta(t)}_{L^p(\Omega)}\leq \norm{\theta(t)}_{L^m(\Omega)}^{m/p}\norm{\theta(t)}_{L^\infty(\Omega)}^{1-m/p}\leq C[\nu(t-s)]^{-(1/m-1/p)}\norm{\theta(s)}_{L^m(\Omega)}.
    \end{align*}
     It remains to remove the assumption $\theta(s)\in L^\infty(\Omega)$. For general $\theta(s)\in L^m(\Omega)$, approximate $\theta(s)$ in $L^m(\Omega)$ by bounded initial data, and let $\theta_n$ be the corresponding solutions. By linearity and Lemma \ref{lp-contraction}, 
    \begin{equation*}
        \sup_{r\in[s,t]}\norm{\theta_n(r)-\theta(r)}_{L^m(\Omega)}\leq\norm{\theta_n(s)-\theta(s)}_{L^m(\Omega)}\to0.
    \end{equation*}
    The assertion follows by lower semicontinuity of the $L^p(\Omega)$ norm.
\end{proof}
\section{Parabolic smoothing for Navier friction}
In this section, we show that parabolic smoothing also holds for the vorticity $\omega^\nu$ with Navier boundary condition, given that the viscosity $\nu$ is sufficiently small. We start with some straightforward observations for the Biot--Savart operator.
\begin{lem}\label{biot-savart-tangential-estimates}
    Let $f\in L^{\infty}(\Omega)$. Then, we have
    \begin{equation*}
        \norm{K_{\Omega}(f)\cdot\tau}_{L^m(\partial\Omega)}\leq C\norm{f}_{L^1(\Omega)}^{1/2}\norm{f}_{L^m(\Omega)}^{1/2} \text{ for all }m\in[1,\infty),
    \end{equation*}
    and
    \begin{equation*}
        \norm{K_{\Omega}(f)\cdot\tau}_{L^\infty(\partial\Omega)}\leq C(m)\norm{f}_{L^m(\Omega)}^{m/2}\norm{f}_{L^\infty(\Omega)}^{1-m/2}\text{ for all }m\in[1,2).
    \end{equation*}
\end{lem}
\begin{proof}
    Without loss of generality, we will assume $f\not\equiv0$. Let $P_{\Omega}:=\partial_nG_{\Omega}$ be the Dirichlet Poisson kernel (for more details, we refer to \cite[Section 1.3]{krantz2}). For $y\in\partial\Omega$, we have
    \begin{equation*}
        K_{\Omega}(f)(y)\cdot\tau(y)=\int_{\Omega}\nabla^{\perp}_yG_{\Omega}(x,y)\cdot \tau(y) f(x)\;dx=\int_{\Omega}\partial_nG_{\Omega}(x,y)f(x)\;dx=\int_{\Omega}P_{\Omega}(x,y)f(x)\;dx.
    \end{equation*}
    From \cite[Theorem 1]{krantz}, we know that
    \begin{equation}\label{poisson-kernel-estimate}
        0\leq P_{\Omega}(x,y)\leq C\frac{d(x)}{\abs{x-y}^2}\leq \frac{C}{\abs{x-y}},
    \end{equation}
    which implies, for every $R>0$,
    \begin{equation}\label{inner-poisson-kernel-estimate}
        \sup_{y\in\partial\Omega}\int_{\Omega\cap B_R(y)}P_{\Omega}(x,y)\;dx\leq C\int_{B_R(0)}\frac{1}{\abs{z}}\;dz\leq CR.
    \end{equation}
    Write $K_{\Omega}(f)\cdot\tau = N_Rf+F_Rf$ with
    \begin{equation*}
        (N_Rf)(y):=\int_{\Omega\cap B_R(y)}P_{\Omega}(x,y)f(x)\;dx\quad\text{and}\quad(F_Rf)(y):=\int_{\Omega\setminus B_R(y)}P_{\Omega}(x,y)f(x)\;dx.
    \end{equation*}
    For $m\in(1,\infty)$, weighted Hölder's inequality and Fubini's theorem give
    \begin{equation*}
        \norm{N_Rf}_{L^{m}(\partial\Omega)}^m\leq (CR)^{m-1}\int_{\Omega}\abs{f(x)}^m\underbrace{\int_{\partial\Omega}P_{\Omega}(x,y)\;dS_y}_{=1}\;dx=(CR)^{m-1}\norm{f}_{L^{m}(\Omega)}^m.
    \end{equation*}
    For the far region, using \eqref{poisson-kernel-estimate} gives $P_{\Omega}^m\leq (C/R)^{m-1}P_{\Omega}$, which together with Minkowski's inequality yields
    \begin{equation*}
        \norm{F_Rf}_{L^{m}(\partial\Omega)}\leq \int_{\Omega}\abs{f(x)}\left(\int_{\partial\Omega\setminus B_R(x)}P_{\Omega}(x,y)^m\;dS_y\right)^{1/m}\;dx\leq C R^{1/m-1}\norm{f}_{L^{1}(\Omega)}.
    \end{equation*}
    Choosing $R=\left(\norm{f}_{L^1(\Omega)}/\norm{f}_{L^m(\Omega)}\right)^{m/(2m-2)}$ proves the first estimate for $m>1$. For $m=1$, it follows directly from Fubini's theorem.\\
    For the second estimate, fix $m\in(1,2)$ and $y\in\partial\Omega$. With \eqref{inner-poisson-kernel-estimate}, we directly get
    \begin{equation*}
        \abs{(N_Rf)(y)}\leq \norm{f}_{L^\infty(\Omega)}\int_{\Omega\cap B_R(y)}P_{\Omega}(x,y)\;dx\leq CR\norm{f}_{L^\infty(\Omega)}.
    \end{equation*}
    For the far region, \eqref{poisson-kernel-estimate}, Hölder's inequality and polar coordinates give
    \begin{align*}
        \abs{(F_Rf)(y)}&\leq C\int_{\Omega\setminus B_R(y)}\frac{\abs{f(x)}}{\abs{x-y}}\;dx\leq C\norm{f}_{L^m(\Omega)}\left(2\pi\int_R^\infty r^{1-m/(m-1)}\;dr\right)^{(m-1)/m}\\&\leq C(m)R^{1-2/m}\norm{f}_{L^m(\Omega)}.
    \end{align*}
    As a remark, we want point out that $C(m)$ diverges as $m\to2$. For $m=1$, the same estimate follows directly from \eqref{poisson-kernel-estimate}. We conclude by choosing\\ $R=\left(\norm{f}_{L^m(\Omega)}/\norm{f}_{L^\infty(\Omega)}\right)^{m/2}.$
\end{proof}
The following lemma is the main driver of the article. It establishes parabolic smoothing rates for $\omega^\nu$ which are sensitive to the scaling of friction $\alpha^\nu$ for sufficiently small viscosity $\nu$. The proof is based on a bootstrap argument, which allows us to first use Lemma \ref{advection-diffusion-lemma} on a startup time interval. We then consider sufficiently small $\nu$ such that the bootstrap hypothesis can be continued to the whole interval $[0,T]$ by continuity.
\begin{lem}\label{navier-parabolic-smoothing}
    Let $\omega_{0}\in H^{1}(\Omega)\cap L^{\infty}(\Omega)$ with $u_{0}=K_{\Omega}(\omega_{0})\in W$, and $\omega^{\nu}= \curl u^\nu$ be the vorticity associated with the strong solution $u^\nu$ of \eqref{navier-stokes} for $u^{\nu}(0,\cdot)=u_0$. Let $r\in(1,2)$ and $q\in[0,\frac{r-1}{2r})$. For $\sigma\in(q,\frac{r-1}{2r})$, there exists a $\nu_0>0$, depending on $\norm{\curl u_0}_{L^r(\Omega)}$, such that for every $\nu\in(0,\nu_0)$ and $t\in(0,T]$,
    \begin{equation*}
        \norm{\omega^\nu(t)}_{L^p(\Omega)}\leq C\nu^{-2\max\set{rq-\sigma,0}/p}(\nu t)^{-(1/r-1/p)}\norm{\curl u_0}_{L^r(\Omega)},\quad r\leq p\leq \infty.
    \end{equation*}
    Here $1/\infty=0$, and the constant $C>0$ is independent of the viscosity $\nu$.
\end{lem}

\begin{proof}
    In the following, we will show the estimates
    \begin{equation}\label{goal-estimates}
        \norm{\omega^\nu(t)}_{L^r(\Omega)}^r\leq C\nu^{-2\max\set{rq-\sigma,0}}\norm{\omega^\nu(0)}_{L^r(\Omega)}^r,\quad(\nu t)^{1/r}\norm{\omega^\nu(t)}_{L^\infty(\Omega)}\leq C\norm{\omega^\nu(0)}_{L^r(\Omega)},
    \end{equation}
 which imply the statement by interpolation
    \begin{align*}
        \norm{\omega^\nu(t)}_{L^p(\Omega)}\leq \norm{\omega^\nu(t)}_{L^r(\Omega)}^{r/p}\norm{\omega^\nu(t)}_{L^\infty(\Omega)}^{1-r/p}\leq C\nu^{-2\max\set{rq-\sigma,0}/p}(\nu t)^{-(1/r-1/p)}\norm{\curl u_0}_{L^r(\Omega)}.
    \end{align*}
    For better readability, we will set
    \begin{equation*}
        X(t):=\sup_{s\in(0,t]}(\nu s)^{1/r}\norm{\omega^\nu(s)}_{L^\infty(\Omega)},\quad X(0):=0.
    \end{equation*}
    Observe that $X\in C([0,T])$. For each $\nu>0$, boundedness of $\omega^\nu$ yields $X(t)\to0$ as $t\to0$. At positive times, we have left continuity as a consequence of lower semicontinuity in $L^\infty(\Omega)$, and right continuity due to the maximum principle (see \cite[Corollary 6.26]{lieberman}) and the continuity of the boundary data in $L^\infty(\partial\Omega)$. Choose $\gamma\in(1/r,1)$ such that $\sigma<\frac{\gamma-1/r}{1+\gamma}$, which ensures 
    \begin{equation*}
        \gamma-\frac{1}{r}-(1+\gamma)\sigma>0.
    \end{equation*}
    Define the stopping time
    \begin{equation*}
        \tau_{\nu}:=\sup\set{t\in(0,T]:\nu^{\gamma-\frac{1}{r}-(1+\gamma)\sigma}X(t)\leq 1},
    \end{equation*}
    which is well defined due to $X\in C([0,T])$ and $X(0)=0$. With this stopping condition, we can establish the estimates \eqref{goal-estimates} on $(0,\tau_\nu]$. An important intermediate step in achieving this is to consider the potentially smaller interval $(0,\min\set{\nu,\tau_\nu}]$ on which we will derive a uniform-in-small-viscosity $L^r(\Omega)$ bound for vorticity. After \eqref{goal-estimates} is shown on $(0,\tau_\nu]$, we close the argument by obtaining a contradiction for $\tau_\nu<T$. \\
    We want to use Lemma \ref{advection-diffusion-lemma} on $[0,\tau_\nu]$, for which we need a suitable $B$ that satisfies \eqref{normal-drift-assumption}. Let $\xi>0$ be small enough that $d\in C^2(\Omega_{\xi})$. Moreover, let $y\in \partial\Omega$ be the nearest boundary point for $x\in\Omega_{\xi}$. Recall that $d(x)=\abs{x-y}$ and $\nabla d(x)=-n(y)$. Then, we obtain for $(s,x)\in(0,\tau_\nu)\times\Omega_\xi$, by Morrey's inequality \ref{morrey} and Lemma \ref{curl_estimate}, that
    \begin{equation}\label{normal-velocity-vorticity}
        \begin{aligned}
            \abs{u^\nu(s,x)\cdot\nabla d(x)}&=\abs{(u^\nu(s,x)-u^\nu(s,y))\cdot\nabla d(x)}\leq \abs{u^\nu(s,x)-u^\nu(s,y)}\\&\leq [u^\nu(s)]_{C^{0,\gamma}}\abs{x-y}^\gamma\leq C \norm{u^\nu(s)}_{W^{1,2/(1-\gamma)}(\Omega)}d(x)^\gamma\\
        &\leq C \norm{\omega^\nu(s)}_{L^{2/(1-\gamma)}(\Omega)}d(x)^\gamma\leq C \norm{\omega^\nu(s)}_{L^{\infty}(\Omega)}d(x)^\gamma\\
        &\leq  \underbrace{C(\nu s)^{-1/r}X(\tau_{\nu})}_{=:B(s)}d(x)^\gamma.
        \end{aligned}
    \end{equation}
    By the stopping condition, we have
    \begin{align*}
        \nu^{\gamma/(1+\gamma)}B(s)^{1/(1+\gamma)}\leq C^{1/(1+\gamma)}\left(\nu^{\gamma-\frac{1}{r}-(1+\gamma)\sigma}X(\tau_{\nu})\right)^{1/(1+\gamma)}\nu^\sigma s^{-1/(r+r\gamma)}\leq C\nu^\sigma s^{-1/(r+r\gamma)}.
    \end{align*}
    Therefore, we get with Lemma \ref{advection-diffusion-lemma} for convex $\beta\in C(\R;[0,\infty))$ that
    \begin{equation}\label{convex-vorticity-estimate}
        \begin{aligned}
            \int_{\Omega}\beta(\omega^\nu(t))\;dx\leq&\int_{\Omega}\beta(\omega^\nu(s))\;dx\\&+C\int_s^t\left(\sqrt{\frac{\nu}{t-z}}+C\nu^\sigma z^{-1/(r+r\gamma)}\right)\int_{\partial\Omega}\beta(\omega^\nu(z))\;dS\;dz,\;
            0\leq s<t\leq\tau_\nu.
        \end{aligned}
    \end{equation}
    For the boundary integral, we later want the interpolation estimates of Lemma \ref{biot-savart-tangential-estimates}. Therefore, we start with an $L^1(\Omega)$ bound for vorticity. We apply \eqref{convex-vorticity-estimate} for $\beta=\abs\cdot$ and Lemma \ref{biot-savart-tangential-estimates} for $m=1$ to obtain for $t\in(0,\tau_{\nu}]$
    \begin{align*}
        \norm{\omega^\nu(t)}_{L^1(\Omega)}&\leq \norm{\omega^\nu(0)}_{L^1(\Omega)}+C\int_0^t\left(\sqrt{\frac{\nu}{t-z}}+C\nu^\sigma z^{-1/(r+r\gamma)}\right)\norm{\omega^\nu(z)}_{L^1(\partial\Omega)}\;dz\\
        &\leq \norm{\omega^\nu(0)}_{L^1(\Omega)}+C\left(\nu^{1/2-q}+\nu^{\sigma-q}\right)\sup_{s\in[0,\tau_\nu]}\norm{\omega^\nu(s)}_{L^1(\Omega)},
    \end{align*}
    which yields
    \begin{equation}\label{vorticity-l1-bound}
        \sup_{s\in[0,\tau_\nu]}\norm{\omega^\nu(s)}_{L^1(\Omega)}\leq C\norm{\omega^\nu(0)}_{L^1(\Omega)},
    \end{equation}
    by choosing $\nu_0>0$ small enough that
    \begin{equation*}
        C\left(\nu^{1/2-q}+\nu^{\sigma-q}\right)\leq \frac{1}{2} \text{ for all }\nu\in(0,\nu_0).
    \end{equation*}
    Next, we want to examine \eqref{convex-vorticity-estimate} for $\beta=\abs{\cdot}^r$. Before doing so, we use Lemma \ref{biot-savart-tangential-estimates} and \eqref{vorticity-l1-bound} to estimate the boundary integral
    \begin{align*}
        \norm{\omega^\nu(s)}^r_{L^r(\partial\Omega)}&\leq \norm{2\kappa-\alpha^\nu}_{L^\infty(\partial\Omega)}^r\norm{K_{\Omega}(\omega^\nu(s))\cdot\tau}^r_{L^r(\partial\Omega)}\leq C\nu^{-rq}\norm{\omega^\nu(s)}_{L^1(\Omega)}^{r/2}\norm{\omega^\nu(s)}_{L^r(\Omega)}^{r/2}\\&\leq C\nu^{-rq}\norm{\omega^\nu(0)}_{L^1(\Omega)}^{r/2}\norm{\omega^\nu(s)}_{L^r(\Omega)}^{r/2}.
    \end{align*}
    Importantly, the estimate is sublinear in $\norm{\omega^\nu}_{L^r(\Omega)}^r$, which allows us to absorb the supremum into the left-hand side later with Young's inequality.  Let $t\in(0,\tau_{\nu}]$. Then, for $s\in(0,t]$, we have
    \begin{equation*}
    \begin{aligned}
        &\norm{\omega^\nu(s)}_{L^r(\Omega)}^r\\&\leq \norm{\omega^\nu(0)}_{L^r(\Omega)}^r+C\nu^{-rq}\norm{\omega^\nu(0)}_{L^1(\Omega)}^{r/2}\sup_{s\in[0,t]}\norm{\omega^\nu(s)}_{L^r(\Omega)}^{r/2}\int_0^s\sqrt{\frac{\nu}{s-z}}+C\nu^\sigma z^{-1/(r+r\gamma)}\;dz\\&\leq\norm{\omega^\nu(0)}_{L^r(\Omega)}^r+C\norm{\omega^\nu(0)}_{L^1(\Omega)}^{r/2}\sup_{s\in[0,t]}\norm{\omega^\nu(s)}_{L^r(\Omega)}^{r/2}\left(\nu^{1/2-rq}t^{1/2}+\nu^{\sigma-rq}t^{1-1/(r+r\gamma)}\right)\\
        &\leq \norm{\omega^\nu(0)}_{L^r(\Omega)}^r+C\norm{\omega^\nu(0)}_{L^1(\Omega)}^r\left(\nu^{1/2-rq}t^{1/2}+\nu^{\sigma-rq}t^{1-1/(r+r\gamma)}\right)^2+\frac{1}{2}\sup_{s\in[0,t]}\norm{\omega^\nu(s)}_{L^r(\Omega)}^r.
    \end{aligned}
    \end{equation*}
    Recall that $rq<\frac{1}{2}$ and $\frac{1}{r+r\gamma}<\frac{1}{2}$. With $1+\nu^{2(\sigma-rq)}\leq 2\nu^{-2\max\set{rq-\sigma,0}}$, we establish the first estimate of \eqref{goal-estimates} on $[0,\tau_\nu]$,
    \begin{equation*}
        \sup_{s\in[0,\tau_{\nu}]}\norm{\omega^\nu(s)}_{L^r(\Omega)}^r\leq C\norm{\omega^\nu(0)}_{L^r(\Omega)}^r\nu^{-2\max\set{rq-\sigma,0}}
    \end{equation*}
    As previously mentioned, we need an uniform-in-small-viscosity bound in $L^r(\Omega)$ on the startup interval $[0,\min\set{\nu,\tau_{\nu}}]$, which will be important for the $L^\infty(\Omega)$ estimate. For $t\leq\nu$, we achieve it with the same estimate as before
    \begin{equation}\label{startup-lr-estimate}
        \sup_{s\in[0,\min\set{\nu,\tau_{\nu}}]}\norm{\omega^\nu(s)}_{L^r(\Omega)}^r\leq C\norm{\omega^\nu(0)}_{L^r(\Omega)}^r\left[1+(\nu^{1-rq}+\nu^{\sigma+1-1/(r+r\gamma)-rq})^2\right]\leq C\norm{\omega^\nu(0)}_{L^r(\Omega)}^r.
    \end{equation}
    To prove the $L^\infty(\Omega)$ bound, we use operator splitting. Let $a^\nu$ be a solution to 
    \begin{equation*}
        \left\{
	\begin{aligned}
		&\partial_{t}a^{\nu} + u^{\nu}\cdot\nabla a^{\nu} = \nu\Delta a^{\nu}\text{ in }(s,t)\times\Omega, \\
        &a^{\nu}=0\text{ on }(s,t)\times\partial\Omega,\;a^\nu(s)=\omega^\nu(s) \text{ in }\Omega,
	\end{aligned}\right.
    \end{equation*}
    and $b^\nu$ be a solution to
    \begin{equation*}
        \left\{
	\begin{aligned}
		&\partial_{t}b^{\nu} + u^{\nu}\cdot\nabla b^{\nu} = \nu\Delta b^{\nu}\text{ in }(s,t)\times\Omega, \\
        &b^{\nu}=(2\kappa-\alpha^\nu)u^\nu\cdot\tau\text{ on }(s,t)\times\partial\Omega,\;b^\nu(s)=0 \text{ in }\Omega.
	\end{aligned}\right.
    \end{equation*}
    Then, by linearity and uniqueness for the advection-diffusion problem with drift $u^\nu$, we have  $\omega^\nu=a^\nu+b^\nu$ on $[s,t]$. The advantage of this is that we can control $a^\nu$ with Lemma \ref{parabolic-smoothing} and $b^\nu$ with the maximum principle. For $s\in(0,T]$, we start the operator splitting at $s/2$ to obtain
    \begin{align*}
        \norm{\omega^\nu(s)}_{L^\infty(\Omega)}&\leq  \norm{a^\nu(s)}_{L^\infty(\Omega)}+\norm{b^\nu(s)}_{L^\infty(\Omega)}\\&\leq C\left(\nu\frac{s}{2}\right)^{-1/r}\norm{a^\nu(s/2)}_{L^r(\Omega)}+\sup_{z\in[s/2,s]}\norm{(2\kappa-\alpha^\nu)u^\nu(z)\cdot \tau}_{L^\infty(\partial\Omega)}\\
        &\leq C\left(\nu s\right)^{-1/r}\norm{\omega^\nu(s/2)}_{L^r(\Omega)}+C\nu^{-q}\sup_{z\in[s/2,s]}\norm{K_{\Omega}(\omega^\nu(z))\cdot\tau}_{L^\infty(\partial\Omega)}\\
        &\leq C\left(\nu s\right)^{-1/r}\norm{\omega^\nu(s/2)}_{L^r(\Omega)}+C\nu^{-q}\sup_{z\in[s/2,s]}\norm{\omega^\nu(z)}_{L^r(\Omega)}^{r/2}\norm{\omega^\nu(z)}_{L^\infty(\Omega)}^{1-r/2}.
    \end{align*}
    With Lemma \ref{biot-savart-tangential-estimates}, we get for $t\in(0,\min\set{\nu,\tau_{\nu}}]$,
    \begin{align*}
        X(t)&=\sup_{s\in(0,t]}(\nu s)^{1/r}\norm{\omega^\nu(s)}_{L^\infty(\Omega)}\\&\leq C\sup_{s\in(0,t]}\norm{\omega^\nu(s/2)}_{L^r(\Omega)}+C\nu^{-q}\sup_{s\in(0,t]}(\nu s)^{1/r}\norm{\omega^\nu(s)}_{L^r(\Omega)}^{r/2}\norm{\omega^\nu(s)}_{L^\infty(\Omega)}^{1-r/2}\\
        &\leq C \norm{\omega^\nu(0)}_{L^r(\Omega)}+C\nu^{-q}\norm{\omega^\nu(0)}_{L^r(\Omega)}^{r/2}(\nu t)^{1/2}\left(\sup_{s\in(0,t]}(\nu s)^{1/r}\norm{\omega^\nu(s)}_{L^\infty(\Omega)}\right)^{1-r/2}\\
        &\leq C \norm{\omega^\nu(0)}_{L^r(\Omega)}+C\nu^{1-q}\norm{\omega^\nu(0)}_{L^r(\Omega)}^{r/2}(X(t))^{1-r/2}\\
        &\leq C \norm{\omega^\nu(0)}_{L^r(\Omega)}+C\nu^{2(1-q)/r}\norm{\omega^\nu(0)}_{L^r(\Omega)}+\frac{1}{2}X(t),
    \end{align*}
    which yields
    \begin{equation*}
        X(t) \leq C\norm{\omega^\nu(0)}_{L^r(\Omega)},\quad t\in(0,\min\set{\nu,\tau_{\nu}}].
    \end{equation*}
    If $\tau_\nu\geq\nu$, we start the operator splitting at $\nu/2$. By the maximum principle, Lemma \ref{lp-contraction}, and Lemma \ref{biot-savart-tangential-estimates}, we derive the following bound for the boundary contribution
    \begin{align*}
        &\sup_{s\in[\nu/2,\tau_\nu]}\norm{b^\nu(s)}_{L^\infty(\Omega)}\\&\leq C\nu^{-q}\sup_{s\in[\nu/2,\tau_\nu]}\norm{K_{\Omega}(\omega^\nu(s))\cdot\tau}_{L^\infty(\partial\Omega)}\\&\leq C\nu^{-q}\left(\sup_{s\in[\nu/2,\tau_\nu]}\norm{K_{\Omega}(a^\nu(s))\cdot\tau}_{L^\infty(\partial\Omega)}+\sup_{s\in[\nu/2,\tau_\nu]}\norm{K_{\Omega}(b^\nu(s))\cdot\tau}_{L^\infty(\partial\Omega)}\right)\\
        &\leq C\nu^{-q}\left(\sup_{s\in[\nu/2,\tau_\nu]}\norm{a^\nu(s)}_{L^r(\Omega)}^{r/2}\norm{a^\nu(s)}_{L^\infty(\Omega)}^{1-r/2}+\sup_{s\in[\nu/2,\tau_\nu]}\norm{b^\nu(s)}_{L^1(\Omega)}^{1/2}\norm{b^\nu(s)}_{L^\infty(\Omega)}^{1/2}\right)\\
        &\leq C\nu^{-q}\left(\norm{\omega^\nu(0)}_{L^r(\Omega)}^{r/2}\norm{\omega^\nu(\nu/2)}_{L^\infty(\Omega)}^{1-r/2}+\sup_{s\in[\frac{\nu}{2},\tau_\nu]}\left(\norm{\omega^\nu(s)}_{L^1(\Omega)}+\norm{a^\nu(s)}_{L^1(\Omega)}\right)^{1/2}\norm{b^\nu(s)}_{L^\infty(\Omega)}^{1/2}\right)\\
        &\leq C\nu^{-q}\left(\norm{\omega^\nu(0)}_{L^r(\Omega)}^{r/2}\left(\left(\nu\frac{\nu}{2}\right)^{-1/r}X(\nu/2)\right)^{1-r/2}+\norm{\omega^\nu(0)}_{L^1(\Omega)}^{1/2}\sup_{s\in[\nu/2,\tau_\nu]}\norm{b^\nu(s)}_{L^\infty(\Omega)}^{1/2}\right)\\
        &\leq C\nu^{1-2/r-q}\norm{\omega^\nu(0)}_{L^r(\Omega)}+C\nu^{-q}\norm{\omega^\nu(0)}_{L^1(\Omega)}^{1/2}\left(\sup_{s\in[\nu/2,\tau_\nu]}\norm{b^\nu(s)}_{L^\infty(\Omega)}\right)^{1/2}.
    \end{align*}
    Hence, we have with Young's inequality
    \begin{equation*}
        \sup_{s\in[\nu/2,\tau_\nu]}\norm{b^\nu(s)}_{L^\infty(\Omega)}\leq C\left(\nu^{1-2/r-q}+\nu^{-2q}\right)\norm{\omega^\nu(0)}_{L^r(\Omega)}.
    \end{equation*}
    For $t\in[\nu,\tau_\nu]$, it holds that $t-\nu/2\geq t/2$. Therefore, we obtain
    \begin{align*}
        (\nu t)^{1/r}\norm{\omega^\nu(t)}_{L^\infty(\Omega)}&\leq (\nu t)^{1/r}\left(C\left(\nu \left(t-\frac{\nu}{2}\right)\right)^{-1/r}\norm{\omega^\nu(\nu/2)}_{L^r(\Omega)}+\norm{b^\nu(t)}_{L^\infty(\Omega)}\right)\\
        &\leq C\left[1+\nu^{1-1/r-q}+\nu^{1/r-2q}\right]\norm{\omega^\nu(0)}_{L^r(\Omega)}\\
        &\leq C\norm{\omega^\nu(0)}_{L^r(\Omega)}.
    \end{align*}
    Recall that $q<\frac{r-1}{2r}$ and $rq<\frac{1}{2}$. In total, we have shown
    \begin{equation*}
        X(t)\leq C\norm{\omega^\nu(0)}_{L^r(\Omega)},\quad t\in[0,\tau_\nu],
    \end{equation*}
    and thus the estimates \eqref{goal-estimates} on $(0,\tau_\nu]$. Let us close with a contradiction for $\tau_\nu<T$. If $\tau_\nu<T$, we would have by continuity and monotonicity
    \begin{equation}\label{early-stopping-time}
        \nu^{\gamma-\frac{1}{r}-(1+\gamma)\sigma}X(\tau_{\nu})= 1.
    \end{equation}
    With the estimates \eqref{goal-estimates} on $(0,\tau_\nu]$, we can choose a sufficiently small $\nu_0$ such that
    \begin{equation*}
        C\norm{\curl u_0}_{L^r(\Omega)}\nu_0^{\gamma-\frac{1}{r}-(1+\gamma)\sigma}\leq \frac{1}{2},
    \end{equation*}
    thus contradicting \eqref{early-stopping-time}. Hence, we have $\tau_\nu=T$ which concludes the proof.
\end{proof}
When $q=0$, Lemma \ref{navier-parabolic-smoothing} establishes that the parabolic scaling is the same as for the heat equation.
\begin{cor}\label{navier-parabolic-smooth-q=0}
    Let $\omega_{0}\in H^{1}(\Omega)\cap L^{\infty}(\Omega)$ with $u_{0}=K_{\Omega}(\omega_{0})\in W$, and $\omega^{\nu}= \curl u^\nu$ be the vorticity associated with the strong solution $u^\nu$ of \eqref{navier-stokes} for $u^{\nu}(0,\cdot)=u_0$ with $q=0$. For $r\in(1,\infty)$, there exists a $\nu_0>0$, depending on $\norm{\curl u_0}_{L^r(\Omega)}$, such that for every $\nu\in(0,\nu_0)$ and $t\in(0,T]$,
    \begin{equation*}
        \norm{\omega^\nu(t)}_{L^p(\Omega)}\leq C(\nu t)^{-(1/r-1/p)}\norm{\curl u_0}_{L^r(\Omega)},\quad r\leq p\leq \infty.
    \end{equation*}
    Here $1/\infty=0$, and the constant $C>0$ is independent of the viscosity $\nu$.
\end{cor}
\begin{proof}
    For $r\in (1,2)$, it is a direct application of Lemma \ref{navier-parabolic-smoothing}. Therefore, let $r\in[2,\infty)$ and set 
    \begin{equation*}
        m:=\frac{2r}{r+1}\in(1,2), \text{ so that }\frac{1}{m}-\frac{1}{2}=\frac{1}{2r}.
    \end{equation*} Since $q=0$, we have, with Lemma \ref{biot-savart-tangential-estimates},
    \begin{equation*}
        \norm{\omega^\nu}_{L^r(\partial\Omega)}\leq C\norm{\omega^\nu}_{L^1(\Omega)}^{1/2}\norm{\omega^\nu}_{L^r(\Omega)}^{1/2}\leq C \norm{\omega^\nu}_{L^r(\Omega)}.
    \end{equation*}
    The argument of Lemma \ref{navier-parabolic-smoothing} holds for  $m$ and any $\sigma\in(0,\frac{r-1}{4r})$. Thus, we can use \eqref{convex-vorticity-estimate} with $\beta=\abs{\cdot}^r$ and derive, with a similar argument as before, that
    \begin{equation*}
        \sup_{s\in[0,T]}\norm{\omega^\nu(s)}_{L^r(\Omega)}^r\leq \norm{\omega^\nu(0)}_{L^r(\Omega)}^r+C(\nu^{1/2}+\nu^\sigma)\sup_{s\in[0,T]}\norm{\omega^\nu(s)}_{L^r(\Omega)}^r
    \end{equation*}
    By choosing $\nu_0$ such that $C(\nu^{1/2}+\nu^\sigma)\leq\frac{1}{2}$, we get
    \begin{equation*}
        \sup_{s\in[0,T]}\norm{\omega^\nu(s)}_{L^r(\Omega)}\leq 2\norm{\omega^\nu(0)}_{L^r(\Omega)}.
    \end{equation*}
    Regarding the $L^\infty(\Omega)$ estimate, we have, with Lemma \ref{biot-savart-tangential-estimates} and Lemma \ref{navier-parabolic-smoothing} for $s\in(0,T]$,
    \begin{align*}
        \norm{K_{\Omega}(\omega^{\nu}(s))\cdot\tau}_{L^\infty(\partial\Omega)}&\leq C\norm{\omega^{\nu}(s)}_{L^m(\Omega)}^{m/2}\norm{\omega^{\nu}(s)}_{L^\infty(\Omega)}^{1-m/2}\\
        &\leq C(\nu s)^{-(1/m-1/2)}\norm{\omega^\nu(0)}_{L^m(\Omega)}\\
        &\leq C(\nu s)^{-1/(2r)}\norm{\omega^\nu(0)}_{L^r(\Omega)}.
    \end{align*}
    Next, we start the operator splitting at $t/2$, which yields
    \begin{align*}
        \norm{\omega^\nu(t)}_{L^\infty(\Omega)}&\leq C(\nu t)^{-1/r}\norm{\omega^\nu(t/2)}_{L^r(\Omega)}+C\sup_{s\in[t/2,t]}\norm{K_{\Omega}(\omega^{\nu}(s))\cdot\tau}_{L^\infty(\partial\Omega)}\\
        &\leq C\left((\nu t)^{-1/r}+(\nu t)^{-1/(2r)}\right)\norm{\omega^\nu(0)}_{L^r(\Omega)}\\
        &\leq C(\nu t)^{-1/r}\norm{\omega^\nu(0)}_{L^r(\Omega)},
    \end{align*}
    where we have chosen $\nu_0$ such that $\nu T\leq 1$. We conclude with interpolation
    \begin{equation*}
        \norm{\omega^\nu(t)}_{L^p(\Omega)}\leq \norm{\omega^\nu(t)}_{L^r(\Omega)}^{r/p}\norm{\omega^\nu(t)}_{L^\infty(\Omega)}^{1-r/p}\leq C(\nu t)^{-(1/r-1/p)}\norm{\curl u_0}_{L^r(\Omega)}.
    \end{equation*}
\end{proof}
\begin{rem}
    If we take $r=p$ in Corollary \ref{navier-parabolic-smooth-q=0}, we get
    \begin{equation*}
        \sup_{t\in[0,T]}\norm{\omega^\nu(t)}_{L^p(\Omega)}\leq C\norm{\curl u_0}_{L^p(\Omega)}\text{ for all }p>1.
    \end{equation*}
    This extends the well-known bound for $p>2$ in \cite[Lemma 3]{nussenzveig}.
\end{rem}
An important consequence of Lemma \ref{navier-parabolic-smoothing} is that it allows us to find a suitable $B$ such that \eqref{normal-drift-assumption} is fulfilled for the Navier--Stokes velocity $u^\nu$.
\begin{cor}\label{normal-velocity-bound-cor}
    Let $\omega_{0}\in H^{1}(\Omega)\cap L^{\infty}(\Omega)$ with $u_{0}=K_{\Omega}(\omega_{0})\in W$, and $\omega^{\nu}= \curl u^\nu$ be the vorticity associated with the strong solution $u^\nu$ of \eqref{navier-stokes} for $u^{\nu}(0,\cdot)=u_0$. For $r\in(1,2)$, $q\in[0,\frac{r-1}{2r})$, and $\gamma\in(0,1)$, there exists a $\nu_0>0$, dependent on $\norm{\curl u_0}_{L^r(\Omega)}$, such that for all $\nu\in(0,\nu_0)$, it holds that
    \begin{equation*}
        \abs{u^\nu(s,x)\cdot \nabla d(x)}\leq C(\nu s)^{-1/r}\norm{\curl u_0}_{L^r(\Omega)}d(x)^\gamma \text{ for almost every }(s,x)\in(0,T)\times\Omega_{\xi},
    \end{equation*}
    where the collar width $\xi>0$ is small enough that $d\in C^2(\Omega_{\xi})$, and the constant $C>0$ is independent of the viscosity $\nu$.
\end{cor}
\begin{proof}
    Choose a small enough $\xi>0$ such that \eqref{normal-velocity-vorticity} in Lemma \ref{navier-parabolic-smoothing} is satisfied. Then, we have with Lemma \ref{navier-parabolic-smoothing}
    \begin{equation*}
        \abs{u^\nu(s,x)\cdot \nabla d(x)}\leq C \norm{\omega^\nu(s)}_{L^\infty(\Omega)}d(x)^\gamma\leq C(\nu s)^{-1/r}\norm{\curl u_0}_{L^r(\Omega)}d(x)^\gamma.
    \end{equation*}
\end{proof}
\section{Uniform vorticity bound}
\begin{prop}\label{compatible-uniform-vorticity-bound}
    Let $\omega_{0}\in H^{1}(\Omega)\cap L^{\infty}(\Omega)$ with $u_{0}=K_{\Omega}(\omega_{0})\in W$, and $\omega^{\nu}= \curl u^\nu$ be the vorticity associated with the strong solution $u^\nu$ of \eqref{navier-stokes} for $u^{\nu}(0,\cdot)=u_0$. For $r\in(1,2)$, $\rho\in(1,\infty)$, and $q\geq0$ satisfying $\rho q<\frac{r-1}{2r}$, there exists a $\nu_0>0$, dependent on $\norm{\curl u_0}_{L^r(\Omega)}$, such that for all $\nu\in(0,\nu_0)$, it holds that
    \begin{equation*}
        \sup_{t\in[0,T]}\norm{\omega^\nu(t)}_{L^\rho(\Omega)}\leq C\norm{\omega^\nu(0)}_{L^\rho(\Omega)},
    \end{equation*}
    where the constant $C>0$ is independent of the viscosity $\nu$.
\end{prop}
\begin{proof}
    The idea is to use Lemma \ref{advection-diffusion-lemma}, which is justified due to Corollary \ref{normal-velocity-bound-cor}, and pick a small enough $\nu_0$ such that we can absorb the boundary integral onto the left-hand side. By applying Corollary \ref{normal-velocity-bound-cor}, we have a dependence of $\nu_0$ on $\norm{\curl u_0}_{L^r(\Omega)}$. Since $\rho q<\frac{r-1}{2r}$, we can choose a $\gamma\in(0,1)$ in the following interval
    \begin{equation*}
        \frac{\rho q+1/r}{1-\rho q}<\gamma<1.   
    \end{equation*}
    This ensures
    \begin{equation}\label{gamma-r-inequality}
        \frac{\gamma-1/r}{1+\gamma}-\rho q>0.
    \end{equation}
    With Lemma \ref{advection-diffusion-lemma} for $\beta=\abs{\cdot}^\rho$ and $t\in(0,T]$, we obtain
    \begin{align*}
        \norm{\omega^\nu(t)}_{L^\rho(\Omega)}^\rho \leq \norm{\omega^\nu(0)}^\rho_{L^\rho(\Omega)}+C\int_0^t&\left[\sqrt{\frac{\nu}{t-s}}+\nu^{\gamma/(1+\gamma)}\left(C(\nu s)^{-1/r}\norm{\curl u_0}_{L^r(\Omega)}\right)^{1/(1+\gamma)}\right]\\&\norm{(2\kappa-\alpha^\nu)u^\nu(s)\cdot\tau}^\rho_{L^\rho(\partial\Omega)}\;ds
    \end{align*}
    Regarding the boundary term, we use Lemma \ref{biot-savart-tangential-estimates} to get
    \begin{align*}
        \norm{(2\kappa-\alpha^\nu)u^\nu(s)\cdot\tau}^\rho_{L^\rho(\partial\Omega)}&=\norm{(2\kappa-\alpha^\nu)K_\Omega(\omega^\nu)(s)\cdot\tau}^\rho_{L^\rho(\partial\Omega)}\\&\leq \norm{2\kappa-\alpha^\nu}^\rho_{L^\infty(\partial\Omega)}\norm{K_\Omega(\omega^\nu)(s)\cdot\tau}^\rho_{L^\rho(\partial\Omega)}\\
        &\leq C\nu^{-\rho q}\norm{\omega^\nu(s)}^\rho_{L^\rho(\Omega)},
    \end{align*}
    which yields for the time integral
    \begin{align*}
        &\int_0^t\left[\sqrt{\frac{\nu}{t-s}}+\nu^{\gamma/(1+\gamma)}\left(C(\nu s)^{-1/r}\norm{\curl u_0}_{L^r(\Omega)}\right)^{1/(1+\gamma)}\right]\norm{(2\kappa-\alpha^\nu)u^\nu(s)\cdot\tau}^\rho_{L^\rho(\partial\Omega)}\;ds\\
        &\leq C\nu^{-\rho q}\sup_{t\in[0,T]}\norm{\omega^\nu(t)}^\rho_{L^\rho(\Omega)}\int_0^t\sqrt{\frac{\nu}{t-s}}+\nu^{\gamma/(1+\gamma)}(\nu s)^{-1/(r+r\gamma)}\norm{\curl u_0}_{L^r(\Omega)}^{1/(1+\gamma)}\;ds\\
        &\leq C\sup_{t\in[0,T]}\norm{\omega^\nu(t)}^\rho_{L^\rho(\Omega)}\left(\nu^{1/2-\rho q}t^{1/2}+\nu^{(\gamma-1/r)/(1+\gamma)-\rho q}\norm{\curl u_0}_{L^r(\Omega)}^{1/(1+\gamma)}t^{1-1/(r+r\gamma)}\right)\\
        &\leq C\sup_{t\in[0,T]}\norm{\omega^\nu(t)}^\rho_{L^\rho(\Omega)}\left[\nu^{1/2-\rho q}+\nu^{(\gamma-1/r)/(1+\gamma)-\rho q}\norm{\curl u_0}_{L^r(\Omega)}^{1/(1+\gamma)}\right].
    \end{align*}
    The first exponent of $\nu$ is positive because of $\rho q<\frac{1}{4}$, and the second exponent is positive because of \eqref{gamma-r-inequality}. Thus, there exists a $\nu_0>0$ such that for all $\nu\in(0,\nu_0)$
    \begin{equation}\label{startup-estimate}
        C\left[\nu^{1/2-\rho q}+\nu^{(\gamma-1/r)/(1+\gamma)-\rho q}\norm{\curl u_0}_{L^r(\Omega)}^{1/(1+\gamma)}\right]\leq \frac{1}{2}.
    \end{equation}
    As the right-hand side is independent of $t\in(0,T]$, the estimate holds uniformly in time
    \begin{equation*}
        \sup_{t\in[0,T]}\norm{\omega^\nu(t)}_{L^\rho(\Omega)}^\rho\leq\norm{\omega^\nu(0)}^\rho_{L^\rho(\Omega)}+\frac{1}{2}\sup_{t\in[0,T]}\norm{\omega^\nu(t)}^\rho_{L^\rho(\Omega)},
    \end{equation*}
    which concludes the proof.
\end{proof}
We use an Aubin--Lions argument, as suggested in \cite{nussenzveig}, to transfer the uniform vorticity estimate, the parabolic smoothing estimate, and the normal-drift estimate to weak solutions.
\begin{lem}\label{vorticity-bound-weak}
    Let $\omega_0\in L^p(\Omega)$ with $p\in(1,\infty)$, and $\omega^\nu=\curl u^\nu$ be the vorticity associated with the unique weak solution $u^\nu$ of \eqref{navier-stokes} for $u^\nu(0,\cdot)=K_\Omega(\omega_0)$. For $\rho\in(1,p]$ and $q\geq 0$ satisfying $\rho q<\min\set{\frac{p-1}{2p},\frac{1}{4}}$, there exists a $\nu_0>0$, depending on $\norm{\omega_0}_{L^p(\Omega)}$, such that for all $\nu\in(0,\nu_0)$, it holds that
    \begin{equation*}
        \sup_{t\in[0,T]}\norm{\omega^\nu(t)}_{L^\rho(\Omega)}\leq C\norm{\omega_0}_{L^\rho(\Omega)}.
    \end{equation*}
     Furthermore, for every $\gamma\in(0,1)$ and $r\in(1,\min\set{p,2})$ with $\rho q<\frac{r-1}{2r}$, we have, for all sufficiently small $\nu>0$,
     \begin{equation*}
         \norm{\omega^\nu(t)}_{L^\infty(\Omega)}\leq C(\nu t)^{-1/r}\norm{\omega_0}_{L^r(\Omega)}\text{ for every }t\in(0,T],
     \end{equation*}
     and
     \begin{equation*}
        \abs{u^\nu(s,x)\cdot \nabla d(x)}\leq C(\nu s)^{-1/r}\norm{\omega_0}_{L^r(\Omega)}d(x)^\gamma \text{ for almost every }(s,x)\in(0,T)\times\Omega_{\xi},
    \end{equation*}
    where the collar width $\xi>0$ is small enough that $d\in C^2(\Omega_{\xi})$. The constants $C>0$ are independent of the viscosity $\nu$.
\end{lem}
\begin{proof}
    Let $\omega^\nu_{0,n}\in H^1(\Omega)\cap L^\infty(\Omega)$ be a sequence of compatible initial vorticities such that
    \begin{equation*}
        \omega_{0,n}^\nu\to\omega_0 \text{ in }L^p(\Omega),
    \end{equation*}
    and denote by $u^\nu_n$ the associated strong solution with vorticity $\omega_n^\nu$. As a note, the approximation depends on $\nu$ as the friction coefficient $\alpha^\nu$ depends on it. In order to use Proposition \ref{compatible-uniform-vorticity-bound}, we choose $r\in(1,\min\set{2,p})$ such that $\rho q<\frac{r-1}{2r}$. Then, by the convergence of initial data, there exists a $\nu_0$, depending on $\norm{\omega_0}_{L^p(\Omega)}$, such that
    \begin{equation}\label{n-uniform-bound}
        \sup_{t\in[0,T]}\norm{\omega_n^\nu(t)}_{L^\rho(\Omega)}\leq C\norm{\omega_{0,n}^\nu}_{L^\rho(\Omega)},\quad \nu\in(0,\nu_0),
    \end{equation}
    where the constant $C>0$ is independent of $n$ and $\nu$. In the following, fix $\nu\in(0,\nu_0)$. We have, with \eqref{n-uniform-bound}, Lemma \ref{curl_estimate}, and Poincar\'{e}'s inequality that
    \begin{equation*}
        \sup_n\norm{u^\nu_n}_{L^\infty(0,T;W^{1,\rho}(\Omega))}\leq C\norm{\omega_0}_{L^p(\Omega)}.
    \end{equation*}
    By the argument in the proof of \cite[Proposition 1]{nussenzveig}, with the suggested adaptations, it follows from the energy inequality that
    \begin{equation*}
        \sup_n\norm{\partial_t u^\nu_n}_{L^2(0,T;H^{-1}(\Omega))}\leq C\left(\nu,\norm{K_{\Omega}(\omega_0)}_{L^2(\Omega)}\right)\leq C\left(\nu,\norm{\omega_0}_{L^p(\Omega)}\right).
    \end{equation*}
    We use the Aubin--Lions lemma \ref{aubin-lions} with $X_0=W^{1,\rho}(\Omega)$, $X=L^2(\Omega)$, and $X_1= H^{-1}$ to extract a subsequence $n_k\to\infty$ such that
    \begin{equation*}
        u_{n_k}^\nu\to v \text{ strongly in } C([0,T];L^2(\Omega))\text{ as }k\to\infty,
    \end{equation*}
    for some $v\in C([0,T];L^2(\Omega))$. Without loss of generality, we can assume that the sequence also converges weakly in $L^2(0,T;H^1(\Omega))$ and thus $v\in L^2(0,T;H^1(\Omega))$. We can easily pass to the limit in each term of the weak formulation \eqref{weak-formulation} to conclude that $v=u^\nu$. Furthermore, we have by the strong convergence of velocity that
    \begin{equation*}
        \omega_{n_k}^\nu(t)=\curl u_{n_k}^\nu(t)\to\curl u^\nu(t) \text{ in distributions for }t\in[0,T],
    \end{equation*}
    and hence also weak convergence in $L^\rho(\Omega)$ by \eqref{n-uniform-bound}. With lower semicontinuity, we obtain that
    \begin{equation*}
        \norm{\omega^\nu(t)}_{L^\rho(\Omega)}\leq \liminf_{k\to\infty}\norm{\omega_{n_k}^\nu(t)}_{L^\rho(\Omega)}\leq C\lim_{k\to\infty}\norm{\omega^\nu_{0,n_k}}_{L^\rho(\Omega)}=C\norm{\omega_0}_{L^\rho(\Omega)}.
    \end{equation*}
    After decreasing $\nu_0$ if necessary, Lemma \ref{navier-parabolic-smoothing} yields an $L^\infty(\Omega)$ bound independent of $k$ for each fixed $t\in(0,T]$. Again, with lower semicontinuity, we get
    \begin{align*}
        \norm{\omega^\nu(t)}_{L^\infty(\Omega)}&\leq\liminf_{k\to\infty}\norm{\omega^{\nu}_{n_k}(t)}_{L^\infty(\Omega)}\leq C(\nu t)^{-1/r}\lim_{k\to\infty}\norm{\omega_{0,n_k}^\nu}_{L^r(\Omega)}=C(\nu t)^{-1/r}\norm{\omega_{0}}_{L^r(\Omega)}.
    \end{align*}
    For a further subsequence, which we label with the same index, we also have $u^{\nu}_{n_k}\to u^\nu$ almost everywhere on $(0,T)\times\Omega$. With Corollary 5.5 and the convergence of initial data, we conclude that, for almost every $(s,x)\in (0,T)\times\Omega_\xi$, 
    \begin{align*}
        \abs{u^\nu(s,x)\cdot\nabla d(x)}=\lim_{k\to\infty}\abs{u^\nu_{n_k}(s,x)\cdot \nabla d(x)}&\leq \lim_{k\to\infty}C(\nu s)^{-1/r}\norm{\omega_{0,n_k}^\nu}_{L^r(\Omega)}d(x)^\gamma\\&=C(\nu s)^{-1/r}\norm{\omega_{0}}_{L^r(\Omega)}d(x)^\gamma.
    \end{align*}
\end{proof}
\section{Strong convergence of vorticity}
We are now in a position to show Theorem \ref{main-result}, which we will do in three steps:
\begin{enumerate}
    \item $u^\nu\to u$ strongly in $C([0,T];L^2(\Omega))$;
    
    \item $b^\nu\to 0$ strongly in $C([0,T];L^\rho(\Omega))$;
    \item $a^\nu\to \curl u$ strongly in $C([0,T];L^\rho(\Omega))$.
\end{enumerate}
\begin{proof}[Proof of Theorem \ref{main-result}]
    Throughout the proof, choose $r\in(1,\min\set{p,2})$ such that $\rho q<\frac{r-1}{2r}$.\\
    \textit{Step 1: $u^\nu\to u$ strongly in $C([0,T];L^2(\Omega))$}. By Lemmas \ref{vorticity-bound-weak} and \ref{curl_estimate}, we have 
    \begin{equation*}
        \sup_{\nu\in(0,\nu_0)}\norm{u^\nu}_{L^\infty(0,T;W^{1,\rho}(\Omega))}\leq C\left(\norm{\omega_0}_{L^\rho}\right),
    \end{equation*}
    where the constant $C>0$ is independent of $\nu$ and the threshold $\nu_0$ only depends on $\norm{\omega_0}_{L^p(\Omega)}$. For $\varphi\in V\cap W^{1,\infty}(\Omega)$, the weak formulation \eqref{weak-formulation} yields
    \begin{align*}
        \abs{\int_\Omega \partial_t u^\nu\cdot\varphi\;dx}\leq&\norm{u^\nu}_{L^2(\Omega)}^2\norm{\nabla\varphi}_{L^\infty(\Omega)}+\nu\norm{\nabla u^\nu}_{L^1(\Omega)}\norm{\nabla\varphi}_{L^\infty(\Omega)}\\&+\nu\norm{\kappa-\alpha^\nu}_{L^\infty(\partial\Omega)}\norm{ u^\nu}_{L^1(\partial\Omega)}\norm{\varphi}_{L^\infty(\partial\Omega)}\\
        \leq& C(1+\nu+\nu^{1-q})\norm{\varphi}_{W^{1,\infty}(\Omega)}.
    \end{align*}
    Therefore, Aubin--Lions lemma \ref{aubin-lions} gives, after passing to a subsequence, a limit $u\in C([0,T];H)\cap L^\infty(0,T;W^{1,\rho}(\Omega))$ such that
    \begin{equation*}
        u^{\nu}\to u \text{ strongly in }C([0,T];L^2(\Omega)) \text{ as }\nu\to 0.
    \end{equation*}
    Also, the convergence and uniform bounds are sufficient to pass to the limit in \eqref{weak-formulation}; thus $u$ is a weak solution of \eqref{euler}.\\
    \textit{Step 2: $b^\nu\to 0$ strongly in $C([0,T];L^\rho(\Omega))$}. Let $a^\nu$ solve \eqref{initial-vorticity} and set $b^\nu=\omega^\nu-a^\nu$. With the uniform vorticity bound in Lemma \ref{vorticity-bound-weak}, Lemma \ref{curl_estimate}, and the trace inequality, we have
    \begin{equation}\label{weak-b-estimate}
        \norm{b^\nu}^\rho_{L^\infty(0,T;L^\rho(\partial\Omega))}\leq C\nu^{-\rho q}\norm{u^\nu}^\rho_{L^\infty(0,T;W^{1,\rho}(\Omega))}\leq C\nu^{-\rho q}\norm{\omega_0}^\rho_{L^{\rho}(\Omega)}.
    \end{equation}
    As in Proposition \ref{compatible-uniform-vorticity-bound}, we choose a $\gamma\in\left(\frac{\rho q+1/r}{1-\rho q},1\right)$ such that
    \begin{equation*}
        \frac{\gamma-1/r}{1+\gamma}-\rho q>0.
    \end{equation*}
    Thanks to the normal-drift estimate in Lemma \ref{vorticity-bound-weak}, we can also use Lemma \ref{advection-diffusion-lemma} with $\beta=\abs{\cdot}^\rho$ for $b^\nu$. Following the integration in time in Proposition \ref{compatible-uniform-vorticity-bound} and using \eqref{weak-b-estimate}, we get
    \begin{equation*}
        \sup_{t\in[0,T]}\norm{b^\nu(t)}_{L^\rho(\Omega)}^\rho\leq C\norm{\omega_0}^\rho_{L^\rho(\Omega)}\left(\nu^{1/2-\rho q}T^{1/2}+\nu^{(\gamma-1/r)/(1+\gamma)-\rho q}\norm{\omega_0}_{L^r(\Omega)}^{1/(1+\gamma)}T^{1-1/(r+r\gamma)}\right).
    \end{equation*}
    As both powers of $\nu$ are positive, the boundary contribution vanishes in the limit
    \begin{equation*}
        \sup_{t\in[0,T]}\norm{b^\nu(t)}_{L^\rho(\Omega)}^\rho\to 0 \text{ as }\nu\to0.
    \end{equation*}
    \textit{Step 3: $a^\nu\to \curl u$ strongly in $C([0,T];L^\rho(\Omega))$}. By Step 1, we can use Lemma \ref{linear-stability} to obtain
    \begin{equation*}
        a^\nu\to\omega \text{ strongly in }C([0,T];L^\rho(\Omega)),
    \end{equation*}
    where $\omega\in C([0,T];L^\rho(\Omega))$ is the unique renormalized solution to \eqref{transport-equation} with drift $u$ and initial datum $\omega_0$. It remains to show $\omega=\curl u$. By the Biot--Savart law and Lemma \ref{curl_estimate}, we have
    \begin{align*}
        \norm{u^\nu-K_{\Omega}(\omega)}_{C([0,T];W^{1,\rho}(\Omega))}\leq C\norm{a^\nu-\omega}_{C([0,T];L^{\rho}(\Omega))}+C\norm{b^\nu}_{C([0,T];L^{\rho}(\Omega))}\to 0 \text{ as }\nu\to0.
    \end{align*}
    By uniqueness of limits, this yields $K_{\Omega}(\omega)=u$, which implies $\omega=\curl u$. By approximation and renormalization, we have 
    \begin{equation}\label{renormalization-limit}
        \norm{\omega(t)}_{L^p(\Omega)}=\norm{\omega_0}_{L^p(\Omega)} \text{ for every }t\in[0,T].
    \end{equation}
    Then, continuity in $L^\rho(\Omega)$ implies weak continuity in $L^p(\Omega)$. With \eqref{renormalization-limit}, we can use the Radon--Riesz property to upgrade to strong continuity in $L^p(\Omega)$. By Lemma \ref{curl_estimate}, we conclude that $u\in C([0,T];W^{1,p}(\Omega))$.
\end{proof}
The energy conservation property of the limit is a straightforward consequence of Lemma \ref{vorticity-bound-weak}.
\begin{proof}[Proof of Corollary \ref{energy-conservation}]
    Choose $r\in(1,\min\set{p,2})$ with $\rho q<\frac{r-1}{2r}$. Then, we have with Lemma \ref{vorticity-bound-weak} and Lemma \ref{curl_estimate} that
    \begin{align*}
        \nu\int_0^T\norm{\nabla u^\nu(t)}_{L^2(\Omega)}^2\;dt&\leq C\nu\int_0^T\norm{\omega^\nu(t)}_{L^2(\Omega)}^2\;dt\leq C\nu\int_0^T\norm{\omega^\nu(t)}_{L^1(\Omega)}\norm{\omega^\nu(t)}_{L^\infty(\Omega)}\;dt\\
        &\leq C\nu^{1-1/r}\int_0^Tt^{-1/r}\;dt\to 0 \text{ as }\nu\to0.
    \end{align*}
    Similarly, we obtain with Lemma \ref{vorticity-bound-weak} and Lemma \ref{biot-savart-tangential-estimates} that
    \begin{align*}
        \nu\int_0^T\int_{\partial\Omega}(\abs{\kappa}+\abs{\alpha^\nu})\abs{u^\nu(t)\cdot\tau}^2\;dS\;dt&\leq C \nu^{1-q}\int_0^T\norm{u^\nu(t)\cdot\tau}_{L^2(\partial\Omega)}^2\;dt\\& \leq C \nu^{1-q}\int_0^T\norm{\omega^\nu(t)}_{L^1(\Omega)}\norm{\omega^\nu(t)}_{L^2(\Omega)}\;dt\\
        &\leq C\nu^{1-q-1/(2r)}\int_0^Tt^{-1/(2r)}dt\to0 \text{ as }\nu\to0.
    \end{align*}
    Thus, passing to the limit in \eqref{energy-identity} yields
    \begin{equation*}
        \norm{u(t)}_{L^2(\Omega)}=\norm{K_{\Omega}(\omega_0)}_{L^2(\Omega)}.
    \end{equation*}
\end{proof}
\appendix
\section{Auxiliary results}
\begin{prop}[Morrey's inequality \cite{evans}]\label{morrey}
    Let $\Omega\subset\R^{2}$ be a bounded domain with $C^1$ boundary, and let $\gamma\in(0,1)$. There exists a constant $C>0$ such that every $v\in W^{1,2/(1-\gamma)}(\Omega)$ has a representative in $C^{0,\gamma}(\overline{\Omega})$ satisfying
    \begin{equation*}
        [v]_{C^{0,\gamma}(\overline{\Omega})}\leq C\norm{v}_{W^{1,2/(1-\gamma)}(\Omega)}.
    \end{equation*}
\end{prop}
\begin{lem}[Biot--Savart estimate]\label{curl_estimate}
    Let $\Omega\subset\R^{2}$ be a bounded, simply connected domain with a boundary of class $C^2$ and $p\in(1,\infty)$. Then 
    \begin{equation*}
        K_{\Omega}:L^p(\Omega)\to W^{1,p}(\Omega)\cap H
    \end{equation*}
    is a bounded linear operator satisfying, for every $f\in L^p(\Omega)$,
    \begin{equation*}
        \curl K_{\Omega}(f)=f,\quad \norm{K_{\Omega}(f)}_{W^{1,p}(\Omega)}\leq C\norm{f}_{L^{p}(\Omega)}. 
    \end{equation*}
\end{lem}
\begin{proof}
    Set $\psi(y):=\int_{\Omega}G_{\Omega}(x,y)f(x)\;dx$. By \cite[Theorem 9.15 and Lemma 9.17]{gilbarg}, we have 
    \begin{equation*}
        \psi\in W^{2,p}(\Omega)\cap W^{1,p}_0(\Omega),\quad \Delta\psi=f,    
    \end{equation*}
    and 
    \begin{equation*}
        \norm{K_{\Omega}(f)}_{W^{1,p}(\Omega)}=\norm{\nabla^\perp\psi}_{W^{1,p}(\Omega)}\leq C\norm{\psi}_{W^{2,p}(\Omega)}\leq C\norm{f}_{L^p(\Omega)}.
    \end{equation*}
    Since $p>1$, Sobolev embedding and the homogeneous Dirichlet condition give $\nabla^\perp\psi=K_{\Omega}(f)\in H$.
\end{proof}
\begin{prop}[Aubin--Lions lemma \cite{chen}]\label{aubin-lions}
    Let $X_{0}\subset X\subset X_{1}$ be Banach spaces such that the embedding $X_{0}\subset X$ is compact and the embedding $X\subset X_1$ is continuous. Then
    \begin{equation*}
        \set{u\in L^{p}(0,T;X_{0}), \partial_{t}u\in L^{q}(0,T;X_{1})}
    \end{equation*}
    embeds compactly into
    \begin{itemize}
        \item $L^{p}(0,T;X)$, for $p<\infty$ and $q\in[1,\infty]$,
        \item $C([0,T];X)$, for $p=\infty$ and $q\in(1,\infty]$.
    \end{itemize}
\end{prop}
The following lemma adapts the linear stability argument in \cite{ciampa} to bounded domains.
\begin{lem}\label{linear-stability}
    Let $\Omega\subset\R^{2}$ be a bounded domain with $C^2$ boundary, $\rho\in(1,\infty)$, and let
    \begin{equation*}
        v^\nu\in C([0,T];H)\cap L^2(0,T;V),\quad v\in C([0,T];H)\cap L^\infty(0,T;W^{1,\rho}(\Omega)),
    \end{equation*}
    such that
    \begin{equation*}
        v^\nu\to v \text{ strongly in }C([0,T];L^2(\Omega)) \text{ as }\nu\to 0.
    \end{equation*}
    For $\theta_0\in L^\rho(\Omega)$, let $\theta^\nu\in C([0,T];L^\rho(\Omega))$ be the associated solution to \eqref{advection-diffusion equation} with drift $v^\nu$, boundary $g^\nu=0$, and initial datum $\theta^\nu(0)=\theta_0$. Then, we have 
    \begin{equation*}
        \theta^\nu\to\theta \text{ strongly in }C([0,T];L^\rho(\Omega)),
    \end{equation*}
    where $\theta\in C([0,T];L^\rho(\Omega))$ is the unique renormalized solution of
    \begin{equation}\label{transport-equation}
        \left\{
	\begin{aligned}
		&\partial_t\theta + v\cdot\nabla\theta=0\text{ in }(0,T)\times\Omega, \\
        &\theta(0)=\theta_0 \text{ in }\Omega.
	\end{aligned}\right.
    \end{equation}
\end{lem}
\begin{proof}
    First, we adapt the linear stability argument of \cite[Lemma 3.3]{ciampa} for $\theta_0\in L^\infty(\Omega)$ and then use the approximation scheme in \cite[proof of Theorem 3.2, Step 3]{ciampa} to generalize to $\theta_0\in L^\rho(\Omega)$.\\
    Suppose that $\theta_0\in L^\infty(\Omega)$. With Lemma \ref{lp-contraction}, we have
    \begin{equation*}
        \norm{\theta^\nu}_{L^\infty((0,T)\times\Omega)}\leq \norm{\theta_0}_{L^\infty(\Omega)},\quad \sup_{s\in[0,T]}\norm{\theta^\nu(s)}_{L^2(\Omega)}\leq \norm{\theta_0}_{L^2(\Omega)}.
    \end{equation*}
    By Banach--Alaoglu's theorem, there exists a subsequence, which we label again as $\nu$, such that
    \begin{equation*}
        \theta^\nu \overset{\ast}{\rightharpoonup} \theta \text{ weakly-$*$ in }L^\infty((0,T)\times\Omega),
    \end{equation*}
    for some $\theta\in L^\infty((0,T)\times\Omega)$. The strong convergence of the drift allows us to pass to the limit in the weak formulation. Thus, $\theta$ is a weak solution to \eqref{transport-equation}. Since $v\cdot n=0$, we can use the renormalization theorem on bounded domains \cite[Theorem 3.1]{boyer} to deduce that
    \begin{equation*}
        \theta\in C([0,T];L^2(\Omega)),\quad \norm{\theta(t)}_{L^2(\Omega)}=\norm{\theta_0}_{L^2(\Omega)}\text{ for all }t\in[0,T].
    \end{equation*}
    Since every weak-$*$ subsequential limit is this same solution, the weak-$*$ convergence holds for the entire family. For $\varphi\in C^\infty_c(\Omega)$, $t\to\int_{\Omega}\theta^\nu(t)\varphi\;dx$ is uniformly Lipschitz by the uniform bounds. Compactness, identification of the limit and density in $L^2(\Omega)$ therefore give
    \begin{equation}\label{test-function-convergence}
        \sup_{s\in[0,T]}\abs{\int_\Omega (\theta^\nu(s)-\theta(s))\varphi\;dx}\to 0 \text{ for all }\varphi\in L^2(\Omega).
    \end{equation}
    Assume the convergence in $C([0,T];L^2(\Omega))$ does not hold. Then, there exist $\eps>0$, a sequence of viscosities $\nu\to 0$, and $t_{\nu}\in[0,T]$ such that
    \begin{equation}\label{contradiction-hypothesis}
        \norm{\theta^\nu(t_\nu)-\theta(t_\nu)}_{L^2(\Omega)}\geq \eps.
    \end{equation}
    From the compactness of $[0,T]$, we can extract a subsequence $t_{\nu}\to t\in[0,T]$ (again labeled with $\nu$), for which we have
    \begin{equation*}
        \theta^\nu(t_\nu)\rightharpoonup \theta(t) \text{ weakly in }L^2(\Omega),
    \end{equation*}
    where we used \eqref{test-function-convergence} and the continuity of $\theta$. We also have convergence of the norms
    \begin{equation*}
        \norm{\theta(t)}_{L^2(\Omega)}\leq \liminf_{\nu\to 0}\norm{\theta^\nu(t_\nu)}_{L^2(\Omega)}\leq \limsup_{\nu\to 0}\norm{\theta^\nu(t_\nu)}_{L^2(\Omega)}\leq \norm{\theta_0}_{L^2(\Omega)}=\norm{\theta(t)}_{L^2(\Omega)},
    \end{equation*}
    and therefore with Radon--Riesz
    \begin{equation*}
        \theta^\nu(t_\nu)\to\theta(t) \text{ strongly in }L^2(\Omega).
    \end{equation*}
    This yields a contradiction to \eqref{contradiction-hypothesis}:
    \begin{equation*}
        0<\eps\leq\norm{\theta^\nu(t_\nu)-\theta(t_\nu)}_{L^2(\Omega)}\leq \norm{\theta^\nu(t_\nu)-\theta(t)}_{L^2(\Omega)}+\norm{\theta(t)-\theta(t_\nu)}_{L^2(\Omega)}\to 0.
    \end{equation*}
    So, we have
    \begin{equation*}
        \theta^\nu\to\theta \text{ strongly in }C([0,T];L^2(\Omega)).
    \end{equation*}
    By the contraction property in $L^\infty(\Omega)$ and interpolation, this also yields
    \begin{equation*}
        \theta^\nu\to\theta \text{ strongly in }C([0,T];L^\rho(\Omega)).
    \end{equation*}
    For $\theta_0\in L^\rho(\Omega)$, choose a sequence $\theta_{0,k}\in L^\infty(\Omega)$ such that $\theta_{0,k}\to\theta_0$ in $L^\rho$ and denote by $\theta^\nu_k,\theta_k$ the corresponding solutions. We conclude with
    \begin{align*}
        \norm{\theta^\nu-\theta}_{C([0,T];L^\rho(\Omega))}&\leq\norm{\theta^\nu-\theta^\nu_k}_{C([0,T];L^\rho(\Omega))}+\norm{\theta^\nu_k-\theta_k}_{C([0,T];L^\rho(\Omega))}+\norm{\theta_k-\theta}_{C([0,T];L^\rho(\Omega))}\\
        &\leq2\norm{\theta_0-\theta_{0,k}}_{L^\rho(\Omega)}+\norm{\theta^\nu_k-\theta_k}_{C([0,T];L^\rho(\Omega))}\to 0,
    \end{align*}
    where we first take $\nu\to0$ and then $k\to\infty$.
\end{proof}
\section{Counterexample to strong vorticity convergence at $q=1/2$}
The following counterexample shows that Theorem \ref{main-result} cannot hold for any $p\in(1,\infty)$ with $q=\frac{1}{2}$.
\begin{lem}\label{counterexample}
    Let $\Omega=B_1(0)\subset\R^2$, $T>0$, and let $u^\nu$ be the unique weak solution to \eqref{navier-stokes} for $u_0(x)=\binom{-x_2}{x_1}$ with $\alpha^\nu=\nu^{-1/2}$. Then
    \begin{equation*}
        u^\nu\to u_0 \text{ strongly in }C([0,T];L^2(\Omega)),
    \end{equation*}
    where $u(t)=u_0$ is a stationary Euler solution with vorticity $\omega=2$. However, for any $\rho\in(1,\infty]$, there exist a constant $C>0$ and viscosity threshold $\nu_0>0$ such that
    \begin{equation*}
        \inf_{t\in[T/2,T]}\norm{\omega^\nu(t)-2}_{L^\rho(\Omega)}\geq C\nu^{-(\rho-1)/(2\rho)} \text{ for } \nu\in(0,\nu_0),
    \end{equation*}
    where $1/\infty=0$.
\end{lem}
\begin{proof}
    The solution is given by $u^\nu(t,x)=h^\nu(t,\abs{x})x^\perp$, where $h^\nu$ solves
    \begin{equation*}
        \left\{
	\begin{aligned}
		&\partial_t h^\nu=\nu\left(\partial_{rr}h^\nu+\frac3r\partial_rh^\nu\right),\quad(t,r)\in(0,T)\times(0,1)\\
		&\partial_rh^\nu(t,1)+\alpha^\nu h^\nu(t,1)=0, \quad t\in(0,T),\\
        &\partial_rh^\nu(t,0)=0,\qquad h^\nu(0,r)=1.
	\end{aligned}\right.
    \end{equation*}
    With the maximum principle, we have $0\leq h^\nu\leq 1$. Therefore, we obtain by testing with $(h^\nu-1)r^3$ and integration by parts that
    \begin{align*}
    \frac{1}{2}\frac{d}{dt}
    \int_0^1\abs{h^\nu-1}^2r^3\;dr
    +\nu\int_0^1\abs{\partial_rh^\nu}^2r^3\; dr
    &=\sqrt\nu\ h^\nu(t,1)\bigl(1-h^\nu(t,1)\bigr)\le\frac{\sqrt\nu}{4},
    \end{align*}
    which implies
    \begin{equation*}
    \sup_{t\in[0,T]}
    \norm{u^\nu(t)-u_0}_{L^2(B_1(0))}^2
    =\sup_{t\in[0,T]}2\pi\int_0^1\abs{h^\nu(t,r)-1}^2r^3\;dr\leq \pi T\sqrt\nu.
    \end{equation*}
    We next show that the tangential boundary traces remain separated. For $r\in(1/2,1)$, we set
    \begin{equation*}
        \overline {h^\nu}(t,r)=1+A_T-A_Te^{t-(1-r)/\sqrt\nu},\qquad A_T=(2e^T-1)^{-1}.
    \end{equation*}
    Choose a $\nu_0>0$ such that $\sqrt\nu T\leq \frac{1}{2}$ for $\nu\in(0,\nu_0)$. Then, we have 
    \begin{align*}
\overline {h^\nu}(0,r)&\geq1=h^\nu(0,r),
\quad r\in[1/2,1],\\
\overline{h^\nu}(t,1/2)&\geq1\geq h^\nu(t,1/2),
\quad t\in[0,T].
    \end{align*}
    Moreover, $\overline {h^\nu}$ is a supersolution on $(0,T)\times(1/2,1)$, since
    \begin{equation*}
        \left(\partial_t-\nu\partial_{rr}-\frac{3\nu}{r}\partial_r\right)\overline{h^\nu}=\frac{3 A_T\sqrt{\nu}}{r}e^{t-(1-r)/\sqrt{\nu}}\geq0,
    \end{equation*}
    and its boundary values at $r=1$ satisfy
    \begin{equation*}
        \partial_r\overline{h^\nu}(t,1)+\alpha^\nu\overline{h^\nu}(t,1)=\alpha^\nu(1+A_T-2A_T e^t)\geq 0.
    \end{equation*}
    With the comparison principle, we therefore get $h^\nu\leq\overline{h^\nu}$ on $[0,T]\times[1/2,1]$, and in particular, for $t\in[T/2,T]$,
    \begin{equation*}
        h^\nu(t,1)\leq \overline{h^\nu}(t,1)=1- A_T(e^{t}-1),
    \end{equation*}
    which implies
    \begin{equation*}
        u(t,x)\cdot\tau(x)-u^\nu(t,x)\cdot\tau(x)=1-h^\nu(t,1)\geq A_T(e^{t}-1)\geq A_T(e^{T/2}-1)>0 \text{ on }\partial\Omega.
    \end{equation*}
    Fix $t\in[T/2,T]$ and set $w(r):=r^2(1-h^\nu(t,r))$. Since $w\geq0$, $w(0)=0$ and $w(1)\geq A_T(e^{T/2}-1)>0$, there is a largest radius $r_\star<1$ such that $w(r_\star)= \frac{1}{2}A_T(e^{T/2}-1)$. By continuity, $w(r)\geq \frac{1}{2}A_T(e^{T/2}-1)$ for $r\in[r_\star,1]$. For the annulus $A=\set{x\in B_1(0):r_\star<\abs{x}<1}$, the energy estimate gives
    \begin{equation*}
        \frac{(A_T(e^{T/2}-1))^2}{4}\abs{A}\leq2\pi\int_0^1w(r)^2r\;dr=2\pi\int_0^1(1-h^\nu(t,r))^2r^5\;dr\leq\int_0^1(1-h^\nu(t,r))^2r^3\;dr\leq \pi T\sqrt{\nu},
    \end{equation*}
    which implies $\abs{A}\leq C \sqrt{\nu}$. Observe that $\omega^\nu=\curl u^\nu=2h^\nu+r\partial_rh^\nu$. Thus, we have $w'(r)=r(2-\omega^\nu(t,r))$. With integration, we therefore get
    \begin{equation*}
        \int_A(2-\omega^\nu(t,x))\;dx=2\pi(w(1)-w(r_\star))\geq \pi A_T(e^{T/2}-1).
    \end{equation*}
    We conclude with Hölder's inequality, for every $\rho\in(1,\infty]$,
    \begin{equation*}
        \norm{\omega^\nu(t)-2}_{L^\rho(\Omega)}\geq \frac{1}{\abs{A}^{1-1/\rho}}\int_A (2-\omega^\nu(t,x))\;dx\geq \frac{\pi A_T(e^{T/2}-1)}{\abs{A}^{1-1/\rho}}\geq C \nu^{-(\rho-1)/(2\rho)}.
    \end{equation*}
\end{proof}


\begin{thebibliography}{10}
\bibitem{boyer}\textsc{F. Boyer}.
\newblock{Trace theorems and spatial continuity properties for the solutions
of the transport equation}, Differential and Integral Equations 18(8), 891--934, 2005.
\bibitem{carlen}E.~A.~\textsc{Carlen}, M.~\textsc{Loss}. \newblock{Optimal smoothing and decay estimates for viscously damped conservation laws, with applications to the 2-D Navier-Stokes equation}, Duke Mathematical Journal 81, 135--157, 1995.
\bibitem{chemetov}N.~V.~\textsc{Chemetov}, F.~\textsc{Cipriano}. \newblock{The inviscid limit for the Navier–Stokes equations with slip condition on permeable walls.} J. Nonlinear Sci., 23(5):731–750, 2013.
\bibitem{chen}
X.~\textsc{Chen}, A.~\textsc{Jüngel}, J.-G.~\textsc{Liu}. \newblock{A note on Aubin-Lions-Dubinski\u{i} lemmas}, Acta Appl. Math., 133:33--43, 2014. 
\bibitem{cheskidov}\textsc{A. Cheskidov, M. C. Lopes Filho, H. J. Nussenzveig Lopes, R. Shvydkoy.} \newblock{Energy conservation in two-dimensional incompressible ideal fluids.} Communications in Mathematical Physics 348, 129--143, 2016.
\bibitem{ciampa}G.~\textsc{Ciampa}, G.~\textsc{Crippa}, S.~\textsc{Spirito}. \newblock{Strong convergence of the vorticity for the 2D Euler
Equations in the inviscid limit.} Arch. Ration. Mech. Anal. 240, no. 1, 295--326, 2021.
\bibitem{chemin}\textsc{J.-Y. Chemin}. \textsc{A remark on the inviscid limit for two-dimensional incompressible fluids}, Communications in Partial Differential Equations 21, 1771--1779, 1996.
\bibitem{clopeau}T.~\textsc{Clopeau}, A.~\textsc{Mikeli\'{c}}, R.~\textsc{Robert}. \newblock{On the vanishing viscosity limit for the 2D incompressible Navier--Stokes equations with the friction type boundary conditions}, Nonlinearity 11,  1625--1636, 1998.
\bibitem{demmel}J.~\textsc{Demmel}, E.~\textsc{Wiedemann}. \newblock{Strong Solutions and Applications to the Inviscid Limit for 2D Navier--Stokes with Navier Slip.} \href{https://arxiv.org/abs/2511.04368}{arXiv:2511.04368}, version 2, 2026.
\bibitem{drivas}\textsc{T. D. Drivas, H. Q. Nguyen}. \newblock{Remarks on the emergence of weak Euler solutions in the vanishing viscosity limit}, Journal of Nonlinear Science 29, 709--721, 2019.
\bibitem{evans}L.~C.~\textsc{Evans}. \newblock{Partial Differential Equations, 2nd ed.} Graduate Studies in Mathematics, vol. 19, American Mathematical Society, Providence, RI, 2010.
\bibitem{nussenzveig}
M.~C.~\textsc{Lopes Filho}, H.~J.~\textsc{Nussenzveig Lopes}, G.~\textsc{Planas}.
\newblock{On the inviscid limit for two-dimensional incompressible flow with Navier friction condition}, SIAM J. Math. Anal., 36, pp. 1130--1141, 2005.
\bibitem{gilbarg}D.~\textsc{Gilbarg}, N.~S.~\textsc{Trudinger}. \newblock{Elliptic Partial Differential Equations of Second Order,} Classics in Mathematics, Springer, reprint of the 1998 edition, 2001.
\bibitem{herrero}\textsc{C. Herrero, G. Tocci, S. Merabia, L. Joly}.
\newblock{Fast increase of nanofluidic slip in supercooled water: the key role of dynamics},
Nanoscale 12, 20396--20403, 2020.
\bibitem{kelliher} J.~P.~\textsc{Kelliher}.
\newblock{Navier--Stokes equations with Navier boundary conditions for a bounded domain in the plane}, SIAM J. Math. Anal., 38, pp. 210--232, 2006. 
\bibitem{kelliher2} J.~P.~\textsc{Kelliher}. \newblock{Observations on the vanishing viscosity limit,} Transactions of the American Mathematical Society 369, 2003--2027, 2017.
\bibitem{kim1}\textsc{N. Kim}. \newblock{Large friction limit and the inviscid limit of 2D Navier--Stokes equations under Navier friction condition}, SIAM Journal on Mathematical Analysis 41, 1653--1663, 2009.
\bibitem{kim2}\textsc{N. Kim}. \newblock{Corrigendum: Large friction limit and the inviscid limit of 2D Navier--Stokes equations under Navier friction condition}, SIAM Journal on Mathematical Analysis 45, 1992--1994, 2013.
\bibitem{krantz}S.~G.~\textsc{Krantz}. \newblock{Calculation and estimation of the Poisson kernel}, Journal of Mathematical Analysis and Applications 302, 143--148, 2005.
\bibitem{krantz2}S.~G.~\textsc{Krantz}. \newblock{Function theory of several complex variables, 2nd edition}, AMS Chelsea Publishing, Providence, RI, 2001.
\bibitem{lieberman}G.~M.~\textsc{Lieberman}. \newblock{Second Order Parabolic Differential Equations,} World Scientific, 1996.
\bibitem{maekawa}\textsc{Y. Maekawa, A. Mazzucato}. \newblock{The inviscid limit and boundary layers for Navier--Stokes flows}, Handbook of Mathematical Analysis in Mechanics of Viscous Fluids, 781–828, 2018.
\bibitem{maekawa2}\textsc{Y. Maekawa, M. Paddick}. \newblock{On convergence criteria for incompressible Navier--Stokes equations with Navier boundary conditions and physical slip rates}, RIMS Kôkyûroku 2038, 9--23, 2017.
\bibitem{mcgrath}\textsc{F. J. McGrath}. \newblock{Nonstationary plane flow of viscous and ideal fluids}, Archive for Rational Mechanics and Analysis 27, 329–348, 1968.
\bibitem{nash}J.~\textsc{Nash}. \newblock{Continuity of solutions of parabolic and elliptic equations,} American Journal of Mathematics 80, 931--954, 1958.
\bibitem{navier}C.~L.~M.~H.~\textsc{Navier}.
\newblock{Mémoire sur les lois du mouvement des fluides}, Mémoires de l’Académie des sciences 6, 389--440, 1827.
\bibitem{nguyen} \textsc{T. T. Nguyen}.
\newblock{The inviscid limit of Navier--Stokes with critical Navier-slip boundary conditions for analytic data,}
SIAM Journal on Mathematical Analysis 52, 4616--4637, 2020.
\bibitem{seis} 
H.~J.~\textsc{Nussenzveig Lopes}, C.~\textsc{Seis}, E.~\textsc{Wiedemann}. \newblock{On the vanishing viscosity limit for 2D
incompressible flows with unbounded vorticity}. Nonlinearity, 34(5):3112--3121, 2021.
\bibitem{paddick}\textsc{M. Paddick}. \newblock{Stability and instability of Navier boundary layers,} Differential and Integral Equations 27, 893--930, 2014.
\bibitem{pelz}\textsc{P. F. Pelz, T. Corneli, S. Mehrnia, M. M. G. Kuhr}.
\newblock{Temperature-dependent wall slip of Newtonian lubricants},
Journal of Fluid Mechanics 948, A8, 2022.
\bibitem{silvestre}L.~\textsc{Silvestre}, V.~\textsc{Vicol}, A.~\textsc{Zlato\v{s}}. \newblock{On the loss of continuity for super-critical drift-diffusion equations}, Archive for Rational Mechanics and Analysis 207, 845--877, 2013.
\end{thebibliography}
\end{document}